\documentclass[reqno,11pt]{amsart}
\usepackage{amsmath}
\usepackage{amsxtra}
\usepackage{amscd}
\usepackage{amsthm}
\usepackage{amsfonts}
\usepackage{amssymb}
\usepackage{eucal}
\usepackage{mathrsfs}
\usepackage[all]{xy}

\let\svthefootnote\thefootnote

\usepackage[dvipsnames]{xcolor}
\newcommand{\red}[1]{{\color{red}#1}}

\usepackage[bookmarks=false]{hyperref}
\usepackage{hyperref}
\usepackage{xr-hyper}

\theoremstyle{plain}
\newtheorem{Thm}[equation]{Theorem}
\newtheorem{Cor}[equation]{Corollary}
\newtheorem{Lem}[equation]{Lemma}
\newtheorem{Prop}[equation]{Proposition}
\newtheorem{Conj}[equation]{Conjecture}

\theoremstyle{definition}
\newtheorem{Def}[equation]{Definition}

\theoremstyle{remark}

\newtheorem{Rem}[equation]{Remark}

\newtheorem{conj}{Conjecture}

\renewcommand{\rm}{\mathrm}

\makeatletter
\renewcommand{\subsection}{\@startsection{subsection}{2}{0pt}{-3ex
plus -1ex minus -0.2ex}{-2mm plus -0pt minus
-2pt}{\normalfont\bfseries}} \makeatother

\numberwithin{equation}{subsection}

\newif\ifShowLabels
\ShowLabelstrue
\newdimen\theight
\def\TeXref#1{%
    \leavevmode\vadjust{\setbox0=\hbox{{\tt
        \quad\quad  {\small \rm #1}}}%
    \theight=\ht0
    \advance\theight by \lineskip
    \kern -\theight \vbox to
    \theight{\rightline{\rlap{\box0}}%
    \vss}%
    }}%

\ShowLabelsfalse

\renewcommand{\sec}[2]{\section{#2}\label{S:#1}%
    \ifShowLabels \TeXref{{S:#1}} \fi}
\newcommand{\ssec}[2]{\subsection{#2}\label{SS:#1}%
    \ifShowLabels \TeXref{{SS:#1}} \fi}

\newcommand{\refs}[1]{Section ~\ref{S:#1}}
\newcommand{\refss}[1]{Section ~\ref{SS:#1}}

\newcommand{\reft}[1]{Theorem ~\ref{T:#1}}
\newcommand{\refl}[1]{Lemma ~\ref{L:#1}}
\newcommand{\refp}[1]{Proposition ~\ref{P:#1}}
\newcommand{\refc}[1]{Corollary ~\ref{C:#1}}

\newcommand{\refe}[1]{\eqref{E:#1}}

\newenvironment{thm}[1]%
    { \begin{Thm} \label{T:#1}  \ifShowLabels \TeXref{T:#1} \fi }%
    { \end{Thm} }

\renewcommand{\th}[1]{\begin{thm}{#1} \sl }
\renewcommand{\eth}{\end{thm} }

\newenvironment{lemma}[1]%
    { \begin{Lem} \label{L:#1}  \ifShowLabels \TeXref{L:#1} \fi }%
    { \end{Lem} }
\newcommand{\lem}[1]{\begin{lemma}{#1} \sl}
\newcommand{\elem}{\end{lemma}}

\newenvironment{propos}[1]%
    { \begin{Prop} \label{P:#1}  \ifShowLabels \TeXref{P:#1} \fi }%
    { \end{Prop} }
\newcommand{\prop}[1]{\begin{propos}{#1}\sl }
\newcommand{\eprop}{\end{propos}}

\newenvironment{corol}[1]%
    { \begin{Cor} \label{C:#1}  \ifShowLabels \TeXref{C:#1} \fi }%
    { \end{Cor} }
\newcommand{\cor}[1]{\begin{corol}{#1} \sl }
\newcommand{\ecor}{\end{corol}}

\newenvironment{defeni}[1]%
    { \begin{Def} \label{D:#1}  \ifShowLabels \TeXref{D:#1} \fi }%
    { \end{Def} }
\newcommand{\defe}[1]{\begin{defeni}{#1} \sl }
\newcommand{\edefe}{\end{defeni}}

\newenvironment{remark}[1]%
    { \begin{Rem} \label{R:#1}  \ifShowLabels \TeXref{R:#1} \fi }%
    { \end{Rem} }
\newcommand{\rem}[1]{\begin{remark}{#1}}
\newcommand{\erem}{\end{remark}}

\newenvironment{conjec}[1]%
    { \begin{Conj} \label{Co:#1}  \ifShowLabels \TeXref{Co:#1} \fi }%
    { \end{Conj} }
\renewcommand{\conj}[1]{\begin{conjec}{#1} \sl }
\newcommand{\econj}{\end{conjec}}

\newcommand{\eq}[1]%
    { \ifShowLabels \TeXref{E:#1} \fi
       \begin{equation} \label{E:#1} }
\newcommand{\eeq}{ \end{equation} }

\newcommand{\prf}{ \begin{proof} }
\newcommand{\epr}{ \end{proof} }

\newcommand\alp{\alpha}     

\newcommand\gam{\gamma}     
\newcommand\del{\delta}     \newcommand\Del{\Delta}
\newcommand\eps{\varepsilon}

\newcommand\tet{\theta}     
\newcommand\iot{\iota}
\newcommand\kap{\kappa}
\newcommand\lam{\lambda}        \newcommand\Lam{\Lambda}

\newcommand\Sit{\Sigma_\tau }
\newcommand\si{\sigma }

\newcommand\calE{{\mathcal{E}}}

\newcommand\calO{{\mathcal{O}}}

\newcommand\calW{{\mathcal{W}}}

\newcommand\C{{\mathbb C}}
\newcommand\sset{\subset}
\newcommand\sminus{\smallsetminus}
\newcommand\en{{\enspace}}

\newcommand{\SC}{\mathcal{C}}

\newcommand{\BD}{\mathbf{D}}

\newcommand\PP{\mathbb{P}}
\renewcommand\AA{\mathbb{A}}

\newcommand\FF{\mathbb{F}}
\newcommand\GG{\mathbb{G}}

\newcommand\ZZ{\mathbb{Z}}

 \newcommand\grg{{\mathfrak{g}}}

\newcommand\sdp{\times \hskip -0.3em {\raise 0.3ex
\hbox{$\scriptscriptstyle |$}}} 

\newcommand\Coker{\operatorname{Coker}}

\newcommand\End{\operatorname{End\,}}
\newcommand\Ext{\operatorname{Ext}}

\newcommand\GL{\operatorname{GL}}

\newcommand\Hom{\operatorname {Hom}}

\newcommand\id{\operatorname{id}}

\newcommand\Int{\operatorname{Int}}

\newcommand\Ker{\operatorname{Ker}}
\newcommand\PGL{{\mathrm{PGL}}}

\newcommand\Proj{\operatorname{Proj}}

\newcommand\rk{\operatorname{rk}}

\newcommand\SL{{\rm SL}}

\newcommand\Spec{\operatorname{Spec}}

\newcommand\supp{\operatorname{supp}}

\newcommand\Tr{\operatorname{Tr}}

\newcommand{\into}{\,\hookrightarrow\,}

\newcommand{\mto}{\mapsto}
\newcommand{\onto}{\,\,\twoheadrightarrow\,\,}

\newcommand\x{\times}
\newcommand\ten{\otimes}

\newcommand{\ra}{\rangle}
\newcommand{\la}{\langle}

\renewcommand\Spec{\operatorname{Spec}}

\newcommand\nc{\newcommand}

\newcommand{\IC}{{\operatorname{IC}}}

\newcommand{\iso}{{\stackrel{\sim}{\longrightarrow}}}

\nc\aff{\operatorname{aff}}
\nc\oGr{\overline{\Gr}}
\nc\Bun{\operatorname{Bun}}
\nc\hgrg{\widehat{\grg}}
\renewcommand\Int{\operatorname{Int}}
\nc\bInt{\overline{\Int}}
\nc\hatLam{\widehat{\Lam}}
\nc\bmu{\overline{\mu}}
\nc\bnu{\overline{\nu}}
\nc\blambda{\overline{\lam}}
\nc\btau{{{t}}}
\nc\bseta{{\boldsymbol{\eta}}}
\renewcommand\SL{\operatorname{SL}}
\nc\ocalW{\overline{\calW}}
\nc\pos{\operatorname{pos}}
\nc\IH{\operatorname{IH}}
\nc\Rep{\operatorname{Rep}}
\nc\Gal{\operatorname{Gal}}
\nc{\tilGr}{\widetilde{\Gr}}

\nc\Pic{\operatorname{Pic}}

\nc{\HC}{{\mathcal{HC}}}
\nc{\on}{\operatorname}
\nc{\BA}{{\mathbb{A}}}
\nc{\BC}{{\mathbb{C}}}
\nc{\BG}{{\mathbb{G}}}
\nc{\BM}{{\mathbb{M}}}
\nc{\BMt}{\BM_\tau}
\nc{\BN}{{\mathbb{N}}}
\nc{\BQ}{{\mathbb{Q}}}
\nc{\BP}{{\mathbb{P}}}
\nc{\BR}{{\mathbb{R}}}
\nc{\BZ}{{\mathbb{Z}}}
\nc{\BS}{{\mathbb{S}}}

\nc{\CA}{{\mathcal{A}}}
\nc{\CB}{{\mathcal{B}}}
\nc{\CalC}{{\mathcal C}}
\nc{\CalD}{{\mathcal D}}
\nc{\CE}{{\mathcal{E}}}
\nc{\CF}{{\mathcal{F}}}
\nc{\CG}{{\mathcal{G}}}
\nc{\CH}{{\mathcal{H}}}
\nc{\CK}{{\mathcal{K}}}
\nc{\CL}{{\mathcal{L}}}
\nc{\CM}{{\mathcal{M}}}
\nc{\CMM}{{\mathcal{M}^{\operatorname{gen}}_\hbar(-\rho)}}
\nc{\CN}{{\mathcal{N}}}
\nc{\CO}{{\mathcal{O}}}
\nc{\CP}{{\mathcal{P}}}
\nc{\CQ}{{\mathcal{Q}}}
\nc{\CR}{{\mathcal{R}}}
\nc{\CS}{{\mathcal{S}}}
\nc{\CT}{{\mathcal{T}}}
\nc{\CU}{{\mathcal{U}}}
\nc{\CV}{{\mathcal{V}}}
\nc{\CW}{{\mathcal{W}}}
\nc{\CX}{{\mathcal{X}}}
\nc{\CY}{{\mathcal{Y}}}
\nc{\CZ}{{\mathcal{Z}}}

\nc{\gen}{{\operatorname{gen}}}
\nc{\cM}{{\check{\mathcal M}}{}}
\nc{\csM}{{\check{\mathcal A}}{}}
\nc{\obM}{{^G\!{\mathsf M}}}
\nc{\obMt}{\obM_\tau}
\nc{\oCA}{{\overset{\circ}{\mathcal A}}{}}
\nc{\obA}{{\overset{\circ}{\mathbf A}}{}}
\nc{\ooM}{{\overset{\circ}{M}}{}}
\nc{\GM}{{^G\!{\mathsf M}}}
\nc{\osR}{{^G\!{\mathsf R}}}
\nc{\GMt}{\GM_\tau}
\nc{\osRt}{\osR_\tau}
\nc{\vM}{{\overset{\bullet}{\mathcal M}}{}}
\nc{\nM}{{\underset{\bullet}{\mathcal M}}{}}
\nc{\obD}{{\overset{\circ}{\mathbf D}}{}}
\nc{\cp}{{\overset{\circ}{\mathbf p}}{}}
\nc{\ofZ}{{\overset{\circ}{\mathfrak Z}}{}}

\nc{\fa}{{\mathfrak{a}}}
\nc{\fb}{{\mathfrak{b}}}
\nc{\fg}{{\mathfrak{g}}}
\nc{\fgl}{{\mathfrak{gl}}}
\nc{\fh}{{\mathfrak{h}}}
\nc{\fri}{{\mathfrak{i}}}
\nc{\fj}{{\mathfrak{j}}}
\nc{\fm}{{\mathfrak{m}}}
\nc{\fn}{{\mathfrak{n}}}
\nc{\fu}{{\mathfrak{u}}}
\nc{\fp}{{\mathfrak{p}}}
\nc{\frr}{{\mathfrak{r}}}
\nc{\fs}{{\mathfrak{s}}}
\nc{\ft}{{\mathfrak{t}}}
\nc{\fT}{{\mathfrak{T}}}
\nc{\ofT}{{\overline{\mathfrak T}}}
\nc{\ofS}{{\overline{\mathfrak S}}}
\nc{\fsl}{{\mathfrak{sl}}}
\nc{\hsl}{{\widehat{\mathfrak{sl}}}}
\nc{\hgl}{{\widehat{\mathfrak{gl}}}}
\nc{\hg}{{\widehat{\mathfrak{g}}}}
\nc{\chg}{{\widehat{\mathfrak{g}}}{}^\vee}
\nc{\hn}{{\widehat{\mathfrak{n}}}}
\nc{\chn}{{\widehat{\mathfrak{n}}}{}^\vee}

\nc{\fA}{{\mathfrak{A}}}
\nc{\fB}{{\mathfrak{B}}}
\nc{\fD}{{\mathfrak{D}}}
\nc{\fE}{{\mathfrak{E}}}
\nc{\fF}{{\mathfrak{F}}}
\nc{\fG}{{\mathfrak{G}}}
\nc{\fI}{{\mathfrak{I}}}
\nc{\fJ}{{\mathfrak{J}}}
\nc{\fK}{{\mathfrak{K}}}
\nc{\fL}{{\mathfrak{L}}}
\nc{\fM}{{\mathfrak{M}}}
\nc{\fN}{{\mathfrak{N}}}
\nc{\frP}{{\mathfrak{P}}}
\nc{\fS}{{\mathfrak S}}
\nc{\fU}{{\mathfrak{U}}}
\nc{\fZ}{{\mathfrak{Z}}}

\nc{\bb}{{\mathbf{b}}}
\nc{\bc}{{\mathbf{c}}}
\nc{\be}{{\mathbf{e}}}
\nc{\bj}{{\mathbf{j}}}
\nc{\bn}{{\mathbf{n}}}
\nc{\bp}{{\mathbf{p}}}
\nc{\bq}{{\mathbf{q}}}
\nc{\bv}{{\mathbf{v}}}
\nc{\bx}{{\mathbf{x}}}
\nc{\by}{{\mathbf{y}}}
\nc{\bw}{{\mathbf{w}}}
\nc{\bA}{{\mathbf{A}}}
\nc{\bB}{{\mathbf{B}}}
\nc{\bC}{{\mathbf{C}}}
\nc{\bK}{{\mathbf{K}}}
\nc{\bD}{{\mathbf{D}}}
\nc{\bH}{{\mathbf{H}}}
\nc{\bM}{{\mathbf{M}}}
\nc{\bN}{{\mathbf{N}}}
\nc{\bS}{{\mathbf{S}}}
\nc{\bT}{{\mathbf{T}}}
\nc{\bV}{{\mathbf{V}}}
\nc{\bW}{{\mathbf{W}}}
\nc{\bX}{{\mathbf{X}}}
\nc{\bP}{{\mathbf{P}}}
\nc{\bQ}{{\mathbf{Q}}}
\nc{\bZ}{{\mathbf{Z}}}

\nc{\sA}{{\mathsf{A}}}
\nc{\sB}{{\mathsf{B}}}
\nc{\sC}{{\mathsf{C}}}
\nc{\sD}{{\mathsf{D}}}
\nc{\sF}{{\mathsf{F}}}
\nc{\sK}{{\mathsf{K}}}
\nc{\sM}{{\mathsf{M}}}
\nc{\sO}{{\mathsf{O}}}
\nc{\sQ}{{\mathsf{Q}}}
\nc{\sP}{{\mathsf{P}}}
\nc{\sV}{{\mathsf{V}}}
\nc{\sW}{{\mathsf{W}}}
\nc{\sZ}{{\mathsf{Z}}}
\nc{\sfp}{{\mathsf{p}}}
\nc{\sr}{{\mathsf{r}}}
\nc{\sfb}{{\mathsf{b}}}
\nc{\sfc}{{\mathsf{c}}}
\nc{\sd}{{\mathsf{d}}}
\nc{\sg}{{\mathsf{g}}}
\nc{\sfl}{{\mathsf{l}}}

\nc{\BK}{{\bar{K}}}

\nc{\tA}{{\widetilde{\mathbf{A}}}}
\nc{\tB}{{\widetilde{\mathcal{B}}}}
\nc{\tg}{{\widetilde{\mathfrak{g}}}}
\nc{\tG}{{\widetilde{G}}}
\nc{\TM}{{\widetilde{\mathbb{M}}}{}}
\nc{\tO}{{\widetilde{\mathsf{O}}}{}}
\nc{\tU}{{\widetilde{\mathfrak{U}}}{}}
\nc{\TZ}{{\tilde{Z}}}
\nc{\tZ}{\widetilde{Z}{}}
\nc{\tx}{{\tilde{x}}}
\nc{\tbv}{{\tilde{\bv}}}
\nc{\tfP}{{\widetilde{\mathfrak{P}}}{}}
\nc{\tz}{{\tilde{\zeta}}}
\nc{\tmu}{{\tilde{\mu}}}

\nc{\td}{\ddot{\underline{d}}{}}
\nc{\tzeta}{\widetilde{\zeta}{}}
\nc{\hd}{{\widehat{\underline{d}}}}
\nc{\hG}{{\widehat{G}}}
\nc{\hBP}{\widehat{\mathbb P}{}}
\nc{\hQ}{{\widehat{Q}}}
\nc{\UM}{{^U\!{\mathsf M}}}
\nc{\hsR}{{^U\!{\mathsf R}}}
\nc{\UMt}{\UM_\tau}
\nc{\hsRt}{\hsR_\tau}
\nc{\hfM}{\widehat{\mathfrak M}{}}
\nc{\hbM}{{^U\!{\mathsf M}}}
\nc{\hbMt}{\hbM_\tau}
\nc{\hCP}{\widehat{\mathcal P}{}}
\nc{\hCR}{\widehat{\mathcal R}{}}
\nc{\hCS}{{\widehat{\mathcal S}}}
\nc{\hfZ}{\widehat{\mathfrak Z}{}}

\nc{\urho}{\underline{\rho}}
\nc{\uB}{\underline{B}}
\nc{\uC}{{\underline{\mathbb{C}}}}
\nc{\ui}{\underline{i}}
\nc{\ofP}{{\overline{\mathfrak{P}}}}

\nc{\hrho}{{\hat{\rho}}}

\nc{\unl}{\underline}
\nc{\ol}{\overline}
\nc{\one}{{\mathbf{1}}}
\nc{\two}{{\mathbf{t}}}

\nc{\Tot}{{\mathop{\operatorname{\rm Tot}}}}
\nc{\Hilb}{{\mathop{\operatorname{\rm Hilb}}}}
\nc{\CHom}{{\mathop{\operatorname{{\mathcal{H}}\it om}}}}
\nc{\defi}{{\mathop{\operatorname{\rm def}}}}
\nc{\length}{{\mathop{\operatorname{\rm length}}}}

\nc{\Cliff}{{\mathsf{Cliff}}}
\nc{\Fl}{{\mathsf{Fl}}}
\nc{\Fib}{{\mathsf{Fib}}}
\nc{\Coh}{{\mathsf{Coh}}}
\nc{\FCoh}{{\mathsf{FCoh}}}

\nc{\cplus}{{\mathbf{C}_+}}
\nc{\cminus}{{\mathbf{C}_-}}
\nc{\cthree}{{\mathbf{C}_*}}
\nc{\Qbar}{{\bar{Q}}}

\nc{\bh}{{\bar{h}}}
\nc{\bOmega}{{\overline{\Omega}}}
\nc\tGr{\widetilde{\Gr}}

\nc{\seq}[1]{\stackrel{#1}{\sim}}
\nc\ogu{\overline{G/U}}
\nc\chlam{\check{\lam}}

\nc\St{\operatorname{St}}

\nc\uS{\underline{S}}
\nc\QM{\mathcal{QM}}
\nc\BPt{\BP^2_\tau}
\nc\bpt{{\mathbf{Q}}_\tau}
\nc\mm{\mu_M}
\nc\mg{\mu_G}
\nc\mt{\mu_\tet}
\nc\ttt{(\tet,\tet')}
\nc\Mtt{\CM_\tau^{\ttt}}
\nc\tnt{(\tet^0,\tet^1)}
\nc\Mt{\CM_\tau^{\tnt}}
\nc\VE{V_\bullet(E)}
\nc\bu{\bullet}
\nc\WF{W_\bu(F)}
\nc\sk{\enskip}
\nc\dmn{(n-d(d-1)/2, 2n-d^2+r, n-d(d+1)/2)}
\nc\Fix{\GM^T}
\nc\Fixt{\GMt^T}
\newcommand\ic{\operatorname{IC}}
\nc\cm{\operatorname{CM}}
\nc\oCM{{^D{\mathsf M}}}
\nc\icm{\ic(\oCM^n)}
\nc\tlam{{\bar\lam}}
\nc\Sch{{\mathsf{Sch}}}
\nc\pg{{p\ccirc\gamma}}

\nc\ccirc{{{}_{\,{}^{^\circ}}}}

\newcommand{\Gr}{{\mathsf {Gr}}}

\newcommand{\Rt}{\mathsf{R}_\tau}
\newcommand{\aMt}{{}^A\mathsf{M}_\tau}
\newcommand{\aRt}{{}^A\mathsf{R}_\tau}

\begin{document}


\begin{flushleft}
{\textbf{Intersection cohomology of the Uhlenbeck\\ compactification 
of the Calogero-Moser space}}
\end{flushleft}

\vskip 5pt

\let\thefootnote\relax\footnote{\hskip -4mm Michael Finkelberg\\
{\tt fnklberg@gmail.com}\\
Victor Ginzburg\\
{\tt ginzburg@math.uchicago.edu}\\
Andrei Ionov\\
{\tt 8916456@rambler.ru}\\
Alexander
Kuznetsov\\
{\tt akuznet@mi.ras.ru}\\
}

\addtocounter{footnote}{-1}\let\thefootnote\svthefootnote
\begin{flushleft}
{\textbf Michael Finkelberg}
\footnote{National Research University
Higher School of Economics, Russian Federation,\\
Department of Mathematics, 6 Usacheva st., Moscow 119048;\\
Institute for Information Transmission Problems}\ $\bullet$\  {\textbf Victor Ginzburg}\footnote{Department of Mathematics, University of Chicago, Chicago, IL
60637, USA}\ $\bullet$\   {\textbf Andrei Ionov} \footnote{National Research University
Higher School of Economics, Russian Federation,\\
Department of Mathematics, 6 Usacheva st., Moscow119048}\ $\bullet$\  {\textbf Alexander
Kuznetsov}\footnote{Steklov Mathematical Institute, Algebraic Geometry Section,\\
8 Gubkina st., Moscow 119991, Russia;\\
The Poncelet Laboratory, Independent University of Moscow;\\
Laboratory of Algebraic Geometry,\\ 
National Research University Higher School 
of Economics, Russian Federation}
\end{flushleft}
\vskip 12pt

\begin{flushleft}
{\textit{To Joseph Bernstein on his 70th birthday,
with gratitude and admiration}}
\end{flushleft}
\vskip 12pt

\begin{flushleft}
{\textbf{Abstract}}\en\small{
We study the natural Gieseker and Uhlenbeck compactifications of the rational
Calogero-Moser phase space. 
The Gieseker compactification is smooth and provides a small
resolution of the Uhlenbeck compactification. We use the resolution
to  compute
the stalks of the IC-sheaf of the Uhlenbeck compactification.}
\end{flushleft}\vskip 12pt

\begin{flushleft}
{\small{\textsl{I would
say that if one can compute the Poincar\'e polynomial\\
 for  intersection 
cohomology  without a computer  then,\\ probably, there is a small
resolution which gives it.}}\\
(J. Bernstein)}
\end{flushleft}

\vskip 5pt

\sec{int}{Introduction}

\ssec{calmo}{The Calogero-Moser space}
The {\sf Calogero-Moser space} $\sM^n$~\cite{kks} is the quotient modulo a free action
of $\PGL_n$ of the space of pairs of complex $n\times n$-matrices $(X,Y)$ such that
$[X,Y]-\on{Id}$ has rank 1.
The Calogero-Moser space is a smooth connected affine algebraic variety of
dimension $2n$~\cite{W}.



\ssec{uhgi}{The Gieseker and Uhlenbeck compactifications}
More generally, for a parameter $\tau\in\BC^\times$,
 we consider a
graded algebra $A^\tau$ with generators
$x,y,z$, of degree 1, and the following commutation relations 
\begin{equation}\label{A-relations}
[x,z]=[y,z]=0,\quad [x,y]=\tau z^2.
\end {equation}
This algebra is a very special case of the Sklyanin algebras
studied in~\cite{ns}, specifically, it corresponds to the case 
of a degenerate plane cubic curve equal
to a triple line.
We set $\BP^2_\tau=\mathsf{Proj}(A^\tau)$, a non-commutative $\mathsf{Proj}$
in the sense of~\cite{AZ}, see also~\cite{KKO}, and
write $\mathsf{coh}(\BP^2_\tau)=\mathsf{qgr}(A^\tau)$ for the corresponding
abelian category $\mathsf{coh}(\BPt)$ of ``coherent
sheaves''. Associated with an object
$E\in \mathsf{coh}(\BPt)$ there is a well-defined triple
$(r=\rk E,\,d=\deg E,\,n=c_2(E))$,   of nonnegative integers,
the rank, the degree, and the
second Chern class of $E$, respectively.

Given a triple $(r,d,n)$, where $r$ and $d$ are coprime,
we introduce  two different  moduli spaces,  $\GMt(r,d,n)$ and
$\UMt(r,d,n)$, of coherent sheaves on
$\BP^2_\tau$. These moduli spaces are
defined by  stability
conditions. The moduli space
$\GMt(r,d,n)$,
 the {\em Gieseker} moduli space, is defined using   Gieseker stability.
 The moduli space $\UMt(r,d,n)$,  the {\em Uhlenbeck}  moduli space,
 is defined using  Mumford stability.
These moduli spaces are projective varieties
which provide two different compactifications of the
moduli space of locally free sheaves.
The variety $\GMt(r,d,n)$ is a particular
case of  moduli spaces studied in~\cite{ns} (cf. also~\cite{NV})
in greater generality.
The  variety  $\UMt(r,d,n)$ is more mysterious; it does not
fit into the framework of ~\cite{ns} and it has not been considered
there. In fact,  even in the commutative case,
a satisfactory construction of the  Uhlenbeck compactification of 
the moduli space of  locally free sheaves on an arbitrary smooth 
surface
is not known so far, cf. \cite{BFG}.
In the case we are interested in, i.e., in the case of the noncommutative 
surface $\BPt$, the variety  $\UMt(r,d,n)$ will be studied in  Section
2. In particular, using an interpretation of our moduli spaces
in terms of certain moduli spaces of quiver representations,
we construct a projective morphism 
$\gamma_\tau:\ \GMt(r,d,n)\to \UMt(r,d,n)$.
This morphism turns out to be  a resolution of singularities,
provided  $r$ and $d$ are coprime. 

In  this paper, we will mostly be interested 
in the case where   $r=1,\ d=0$, and  $\tau\neq 0$.
The moduli space of locally free sheaves $\CE$ on $\BPt$ such that
$\rk\CE=1,\ \deg\CE=0$, and $c_2(\CE)=n$ has an ADHM type description.
Specifically, according   to~\cite{NS} and \cite{KKO}, this moduli space is
isomorphic to the variety $\sM^n_\tau$ defined as
a  quotient  of the space of pairs  $(X,Y)$,
of $n\times n$-matrices such that
$\rk([X,Y]-\tau\on{Id})=1$, by  the (free) action
of the group $\PGL_n$  by conjugation.
 Note that the rescaling map $(X,Y)\mapsto(\frac{1}{\tau} X,Y)$
gives  a canonical isomorphism of $\sM^n_\tau$
with the Calogero-Moser space
$\sM^n$. Therefore, the varieties $\GM^n_\tau=\GMt(1,0,n)$ and
$\UM^n_\tau=\UMt(1,0,n)$ provide two different  compactifications of the
Calogero-Moser space. Since $1$ and $0$ are coprime,
 the corresponding  morphism $\gamma_\tau:\ \GM^n_\tau\to\UM^n_\tau$ is
 a resolution of singularities.
Moreover, we show that this morphism is {\em small} in the sense
of Goresky--MacPherson.

One can allow the parameter $\tau$ to vary in $\BA^1$.
Similar to the above, one 
constructs the family of Gieseker, resp. Uhlenbeck, compactifications
$\GM^n$,
resp. $\UM^n$,
equipped with maps to $\AA^1$ such that the fibers over the point $\tau \in \AA^1 \setminus \{0\}$
are $\GM^n_\tau$, resp. $\UM^n_\tau$. Furthermore, we construct
a small resolution of singularities $\gamma:\ \GM^n \to \UM^n$. In fact, over $\AA^1 \setminus \{0\}$
the maps $\gamma_\tau:\ \GM^n_\tau\to\UM^n_\tau$ are identified with the maps discussed above, 
while the fiber over $\tau=0$ is
the Hilbert--Chow morphism 
$\gamma_0:\ \on{Hilb}^n\BP^2\to S^n\BP^2=(\BP^2)^n/\fS_n$.

It is well known that $\gamma_0$ is only {\em semismall}. The reason for this
difference between $\gamma_0$ and $\gamma_\tau,\ \tau\ne0$, is due to the
difference between  the stratifications of the commutative and noncommutative
Uhlenbeck compactifications, respectively. Namely, we have a distinguished (classical,
commutative) projective line subscheme $\BP^1\subset\BP^2_\tau$, and
similarly,
we have $\BP^1\subset\BP^2$ such that $\BP^2\setminus\BP^1=\BA^2$.
There is a stratification
\begin{equation*}
\UM^n_\tau=\bigsqcup_{0\leq m\leq n}\sM^m_\tau\times S^{n-m}\BP^1,\qquad\text{for $\tau\ne0$},
\end{equation*}
and $\gamma_\tau$ is an isomorphism over the open part $\sM^n_\tau$.
Similarly, we have a stratification
\begin{equation*}
\UM^n_0=S^n\BP^2=\bigsqcup_{0\leq m\leq n}S^m\BA^2\times S^{n-m}\BP^1,
\end{equation*}
but $\gamma_0$ is {\em not} an isomorphism over $S^n\BA^2$, only a semismall
resolution of singularities.

We remark that the readers experienced with the classical Uhlenbeck
compactifications might have expected a different stratification with strata of  the form
$\sM^m_\tau\times S^{n-m}\BP^2_\tau$
(the reason for semismallness of the classical Gieseker resolution); however $\BP^2_\tau$
is {\em not} a classical scheme, so only its ``classical part'' $\BP^1$
survives in the classical moduli scheme $\UM^n_\tau$, yielding the
stratification of the previous paragraph.

Note also that the Gieseker moduli spaces of~\cite{ns} carry a natural
Poisson structure. We expect it to descend to the Uhlenbeck compactification.
However, even normality of $\UM^n_\tau$ seems to be a hard question and it
is beyond the scope of this paper.


\newcommand{\height}{\operatorname{\mathbf{h}}}

\ssec{gormac}{The main Theorem}


Let $\frP(n)$ denote the set of partitions of an integer $n\geq0$ and for an algebraic variety $T$ and a partition $\lambda = (\lambda_1,\dots,\lambda_l)$ put
\begin{equation*}\label{diagonal-stratification}
S_\lambda T = \{ \textstyle\sum \lambda_i P_i\ |\ P_1 \ne P_2 \ne \dots \ne P_l \in T \} \subset S^nT=T^n/\fS_n,
\end{equation*}
so that 
$S^n T = \bigsqcup_{\lambda \in \frP(n)} S_\lambda T$
is a stratification, which we call the {\sf diagonal stratification}. 
%
Let $\ic(\UM^n_\tau)$ be the IC sheaf (see~\cite{bbd}) 
of the Uhlenbeck compactification.
Our main result is the computation of the stalks of the IC sheaf.


\th{main}
The IC sheaf of the Uhlenbeck compactification is
smooth along the stratifications
\begin{equation*}
\UM^n_\tau = \bigsqcup_{\substack{0 \le m \le n \\[.5ex] \lambda \in \frP(n-m)}}  \left(\sM^m \times S_\lambda\PP^1\right).
\end{equation*}
For $0 \le m \le n$ and $\lambda = (\lambda_1,\dots,\lambda_k) \in \frP(n-m)$,
the stalk of the sheaf $\ic(\UM^n_\tau)$ at a point of a stratum 
$\sM^m\times S_\lambda\BP^1$ is isomorphic to
\begin{equation}\label{heis}
\bigotimes_{i=1}^k \left( \bigoplus_{\mu \in \frP(\lambda_i)} \BC[2l(\mu)] \right)[2m]
\end{equation}
as a graded vector space.
\eth

The proof employs the small resolution of the family
$\gamma:\ \GM^n_\tau\to\UM^n_\tau$ and reduces the study of the fibers
for $\tau\ne0$ to the well-known properties of the fibers of the Hilbert--Chow
morphism for $\tau=0$. 

\begin{Rem} Given a complex semisimple simply connected group $G$, one can 
consider the moduli
space of $G$-bundles on $\BP^2$ equipped with a trivialization
at the infinite line $\BP^1\sset\BP^2$.
There is also an Uhlenbeck space  ${\mathcal U}_G(\BA^2)$ that contains
the above moduli space as a Zariski open subset
(the variety  ${\mathcal U}_G(\BA^2)$ is not proper in this setting).
Assume that the group $G$ is  almost simple,
and let $G_{\on{aff}}$ be the affinization  of $G$, a Kac-Moody group
such that $\fg _{\on{aff}}=\textrm{Lie}(G_{\on{aff}})$ is an affine Lie algebra.
Then, the  Uhlenbeck space  ${\mathcal U}_G(\BA^2)$ may be viewed
as  a slice in the affine Grassmannian for the group  $G_{\on{aff}}$,
that is, in the  double affine Grassmannian for $G$~\cite{bf}.
The IC stalks of ${\mathcal U}_G(\BA^2)$  may be identified, in
accordance with the predictions based on the geometric Satake
correspondence, with certain graded versions
of the weight spaces of the basic integrable representation of
$\fg_{\on{aff}}^\vee$, the Langlands dual of the Lie algebra $\fg _{\on{aff}}$.
In the simply-laced case, the Dynkin diagram of the  Lie algebra
 $\fg _{\on{aff}}$ is an  affine  Dynkin diagram of types
 $\widetilde A,\,\widetilde D,\,\widetilde E$,
and we have $\fg_{\on{aff}}^\vee=\fg_{\on{aff}}$.

It is often useful to view the graph with one vertex and one edge-loop
at that vertex as a  Dynkin diagram  of type $\widetilde A_0$.
It is known that the   Kac-Moody  Lie algebra associated with 
 $\widetilde A_0$ is  the Heisenberg Lie algebra  ${\mathfrak H}$.
By definition, we have ${\mathfrak H} := \C\delta\ltimes
\widehat{\C(\!(t)\!)}$,
where $\widehat{\C(\!(t)\!)}$ is a central extension of the abelian  Lie algebra
$\C(\!(t)\!)$  and $\delta:=t\frac{d}{dt}$, a derivation.
The Fock representation of ${\mathfrak H}$ plays the role of the basic
integrable representation of an affine Lie algebra.
The tensor factors of the graded vector space in \eqref{heis}
may be identified in a natural way with 
certain weight spaces of the Fock space (the action of the derivation
$\delta$ 
gives a grading
on the  Fock space).
This suggests, in view of the above, that
our variety
 $\UM^n_\tau$ might play the role of a slice in
some kind of an affine Grassmannian for the Heisenberg group
and Theorem \ref{T:main} is
a manifestation of (a certain analogue  of) the geometric Satake
correspondence in the case of   Dynkin diagram  of type $\widetilde A_0$.
\end{Rem}

\ssec{orga}{Organization of the paper}
In~\refs{ch5} we study coherent sheaves on a noncommutative projective plane
and the corresponding representations of a Kronecker-type quiver.
We introduce Gieseker and Mumford stabilities of sheaves and interpret them
as stabilities of quiver representations. We construct the Gieseker and
the Uhlenbeck moduli spaces of sheaves as GIT moduli spaces of quiver
representations and a map $\gamma$ between the moduli spaces as the map
coming from a variation of GIT quotients.
In~\refs{ch6} we discuss the special case of sheaves of rank 1 and degree 0.
In this case the Gieseker and the Uhlenbeck moduli spaces are compactifications
of the Calogero-Moser space. We investigate in detail the map $\gamma$ between
the compactifications and compute the stalks of the IC sheaf on the Uhlenbeck compactification.
In the Appendix we provide proofs of some of the results of~\refs{ch5}.
%

\medskip

\noindent
{\textbf{Notation.}}\en Given  a vector space $V$, we write $V^\vee$ for
the dual  vector space and $S^\bullet V=\oplus_{i\geq 0}\ S^i V$  for the
Symmetric algebra of $V$.

\ssec{ackno}{Acknowledgments}
The work of M.F. and A.K. has been funded by the Russian Academic Excellence
Project `5-100'. 
V.G.  was supported in part by the NSF grant DMS-1303462.
A.I. was supported by the grants NSh-5138.2014.1 and RFBR 15-01-09242.
A.K. was partially supported by RFBR 14-01-00416, 15-01-02164, 15-51-50045 and 
by the Simons foundation.





\newcommand{\tH}{{\widetilde{H}}}

\sec{ch5}{Sheaves on the noncommutative plane $\BPt$ and quiver representations}
\noindent
To construct the Gieseker and the Uhlenbeck compactifications of the 
Calogero-Moser space we use an interpretation of the latter as moduli spaces 
of coherent sheaves on a noncommutative projective plane.

\ssec{reco}{Sheaves on the noncommutative projective plane}

\newcommand{\ba}{\mathfrak{a}}

We start with a slightly more invariant definition of the Calogero-Moser space.

We consider a symplectic vector space $H$ of dimension 2 with a symplectic form $\omega \in \Lambda^2H^\vee$,
a vector space $V$ of dimension $n$, a nonzero complex number $\tau \in \C^\times$,
and consider the subvariety
$\widetilde\sM_\tau(V) \subset \Hom(V,V\otimes H)$
defined by
\begin{equation*}
\widetilde\sM_\tau(V) = \{ \ba \in \Hom(V, V\otimes H)\ |\ \mathrm{rank}(\omega( \ba \circ \ba) - \tau\id_V) = 1 \},
\end{equation*}
where $\omega(\ba\circ \ba)$ is defined as the composition
$V \xrightarrow{\ \ba\ } V \otimes H \xrightarrow{\ \ba \otimes \id_H\ } V \otimes H \otimes H \xrightarrow{\ \id_V \otimes \omega\ } V$;
we consider the action of $\PGL(V)$ on $\widetilde\sM_\tau(V)$ by conjugation,
and define
\begin{equation*}
\sM_\tau(V) = \widetilde\sM_\tau(V)/\PGL(V).
\end{equation*}
A choice of  symplectic basis in $H$ allows us to rewrite $\ba$ as a pair of operators $(X,Y)$.
Then $\omega(\ba\circ \ba)$ can be written  as $[X,Y]$,  so this definition
agrees with the standard one.

Let $\tH := H \oplus \BC$. We define a {\sf twisted symmetric algebra of $\tH$} by
\begin{equation*}
A^\tau =
S^\bullet_\tau\tH =
\BC\langle H \oplus \BC z\rangle/\langle [H,z]=0,\ [h_1,h_2]=\tau \omega(h_1,h_2) z^2 \rangle.
\end{equation*}
Choosing a symplectic basis $x,y$ in $H$ the defining relations in
$A^\tau$
take the form \eqref{A-relations}.

%

The algebra $A^\tau$ is a  graded noetherian  algebra 
and we let
\begin{equation*}
\BP^2_\tau := \mathsf{Proj}(A^\tau)
\end{equation*}
be the noncommutative ``projective spectrum'' of 
$A^\tau$ in the sense of~\cite{AZ}. 
The category of ``coherent sheaves'' on the noncommutative scheme
$\BP^2_\tau$ is defined as
$\mathsf{coh}(\BP^2_\tau):=\mathsf{qgr}(A^\tau)$,  a quotient of the
abelian  category
of finitely generated graded $A^\tau$-modules by the 
Serre subcategory of finite-dimensional modules.
Note that the group $\SL(H)$ acts on the algebra $A^\tau$ by automorphisms.
The action on $A^\tau$ induces an $\SL(H)$-action 
on the category of coherent sheaves $\mathsf{coh}(\BPt)$.

As it was shown in~\cite{AZ,KKO,BGK,NS} and other papers, coherent sheaves
on such a noncommutative projective plane behave very similarly to those
on the usual (commutative) plane~$\BP^2$. For instance, one can define
the cohomology spaces of sheaves, local $\underline{\Ext}$ sheaves,
the notions of torsion free and locally free sheaves, one has the sequence
$\{\CO(i)\}_{i\in\ZZ}$ of ``line bundles'', one can prove Serre duality and
construct the Beilinson spectral sequence.

The main differences from the commutative case are
\begin{itemize}
\item
there is {\em no tensor product of sheaves} in general  (due to noncommutativity);
however, one can tensor with $\CO(i)$ and thus define the twist functors $F \mapsto F(i)$
since sheaves $\CO(i)$ correspond to graded $A^\tau$-modules having a natural bimodule structure
(alternatively, the twist functor can be thought of as the twist of the grading functor
in the category of graded $A^\tau$-modules);
\item
the dual of a sheaf on $\PP^2_\tau$ is a sheaf on
$\mathsf{Proj}((A^\tau)^{\mathrm{opp}})$,
the ``opposite'' noncommutative projective plane;
in fact, one has $\mathsf{Proj}((A^\tau)^{\mathrm{opp}})=\PP^2_{-\tau}$ since $(A^\tau)^{\mathrm{opp}} \cong A^{-\tau}$.

\item
the noncommutative projective plane $\BP^2_\tau$ {\em has fewer points} than the usual plane $\BP^2$,
and as a consequence the category $\mathsf{coh}(\BP^2_\tau)$ {\em has more locally free sheaves} than
$\mathsf{coh}(\BP^2)$.
\end{itemize}

%



Below, we summarize the results of \cite{AZ,BGK,KKO,NS} that we are
going to
use later in the paper.

By \cite[Theorem 8.1(3)]{AZ},
the cohomology groups of the sheaves $\CO(i)$ are given by the
following formulas (similar to those in the commutative case):

\begin{equation*}
H^p(\BPt,\CO(i))=
\begin{cases}
A^\tau_i = S^i\tH, &\text{if $p=0$ and $i\ge0$}\\
(A^\tau_{-i-3})^\vee = S^{-i-3}\tH^\vee, &\text{if $p=2$ and $i\le-3$}\\
0 &\text{otherwise.}
\end{cases}
\end{equation*}

%
One has  a functorial Serre duality isomorphism
%
%
\begin{equation*}
\Ext^i(E,F)\cong\Ext^{2-i}(F,E(-3))^\vee.
\end{equation*}
%

The sheaves  $(\CO(-2),\CO(-1),\CO)$ form a full exceptional collection
in the derived category and there is an associated Beilinson type spectral sequence. 
The construction of the  spectral sequence involves the sheaves $Q_0$,
$Q_1$ and $Q_2$ on $\BPt$ defined by
\begin{equation}\label{defq}
Q_0 = \CO,
\qquad
0 \to \CO \xrightarrow{ (x,y,z) } \CO(1) \oplus \CO(1) \oplus \CO(1) \to Q_1 \to 0,
\qquad
Q_2 = \CO(3).
\end{equation}
Sometimes, another resolution for $Q_1$ is more convenient
\begin{equation}\label{anotherq1}
0 \to Q_1 \to \CO(2) \oplus \CO(2) \oplus \CO(2) \xrightarrow{ (x,y,z) } \CO(3) \to 0.
\end{equation}
We remark that each of the two  sequences above is a truncation of the Koszul complex.

The Beilinson spectral sequence 
has the form
\begin{equation*}
E^{-p,q}_1=\Ext^q(Q_{p}(-p),E)\otimes \CO(-p)\Longrightarrow E^i_\infty=\begin{cases} E, &\text{for $i=0$}\\
0, &\text{otherwise}
\end{cases}
\end{equation*}
where $p=0,1,2$. 
Using the Beilinson spectral sequence one shows, cf. \cite[\S7.2]{BGK}, that
any coherent sheaf $E$ on $\BPt$ admits a resolution of the form
\eq{resol1}
0\to V'\otimes \CO(k-2)\to V\otimes \CO(k-1)\to V'' \otimes \CO(k)\to E\to0
\eeq
for some $k\in\BZ$ and  vector spaces $V',V,V''$.

The dual
$
E^* := \underline{\Hom}(E,\CO)$,
of any sheaf $E$ is a  sheaf on the opposite plane $\PP^2_{-\tau}$.
The sheaf $E$ is called {\itshape locally free} if
$\underline{\Ext}^i(E,\CO) = 0$ for $i > 0$. 

The following statements are proved in \cite[Proposition 2.0.4]{BGK}.
For any sheaf $E$, the sheaf $E^*$ is locally free, furthermore,
$E$  is locally free if and only if its canonical map $E\to E^{**}$
is an isomorphism. The kernel of a morphism of locally free sheaves is always locally free.

Let $\CE$ be a  locally free sheaf.
Writing~\refe{resol1} for $\CE^*$ and dualizing, one deduces that
%
%
%
%
any locally free sheaf $\CE$ has a resolution of the form
\eq{resol2}
0\to\CE\to U' \otimes \CO(-k)\to U \otimes \CO(1-k)\to U'' \otimes \CO(2-k)\to0.
\eeq

A sheaf $E$ is called {\itshape torsion free} if it can be embedded in a locally free sheaf.
This can be shown, e.g., using \cite[Proposition
2.0.6]{BGK},  to be equivalent to 
the injectivity of the canonical map $E \to E^{**}$.

For a coherent sheaf $E$ its {\itshape Hilbert polynomial} is defined by
the  usual formula
\begin{equation*}
h_E(t)=\sum_{i=0}^2(-1)^i\dim H^i(\BPt,E(t)).
\end{equation*}
For sheaves $\CO(i)$ it is the same as in the commutative case
$h_{\CO(i)}(t)=(t+i+1)(t+i+2)/2$. So, using~\refe{resol1} one sees
that the Hilbert polynomial of any sheaf can be written as
\eq{david}
h_E(t)=r(E)\frac{(t+1)(t+2)}{2}+\deg(E)\frac{2t+3}{2}+\frac{\deg(E)^2}{2}-c_2(E)
\end{equation}
for some integers $r(E),\deg(E)$ and $c_2(E)$ defined by this equality and
called the rank, degree and second Chern class of $E$ respectively.
It is clear from the definition that the Hilbert polynomial as well as the rank
and the degree are additive in exact sequences. Further, one can check
that they behave naturally with respect to dualization:
\begin{itemize}
\item for any sheaf $E$ one has $r(E^*) = r(E)$;
\item for a torsion free sheaf $E$ one also has $\deg(E^*) = - \deg(E)$;
\item for a locally free sheaf $\CE$ one also has $c_2(E^*) = c_2(E)$.
\end{itemize}
\newcommand{\cch}{\mathop{\mathrm{ch}}\nolimits}

It is more convenient, sometimes, to use the class
\begin{equation*}
\cch_2(E) := \deg(E)^2/2 - c_2(E),
\end{equation*}
(rather than the second Chern class $c_2(E)$), which may be thought of as the second coefficient of the Chern character
and which  is additive in exact sequences.

For any sheaf $E$ the rank $r(E)$ is nonnegative. If $E$ is torsion free and nonzero,
then $r(E) > 0$; moreover,  the degree $\deg(E)$ is nonnegative if $r(E) = 0$.
The sheaf $F$ is called {\itshape Artin sheaf of length} $n = h_F =
\cch_2(F)$
if  both the rank and  degree of  $F$ are equal to zero, equivalently, 
 the Hilbert polynomial
of $F$ is constant. In this case,
 the integer $n := h_F =
\cch_2(F)$ is  nonnegative and it is called  the length of $F$.

A special feature of the noncommutative plane $\BPt$
is that it has fewer points than the commutative $\BP^2$: all points of
 $\BPt$
are contained, in a sense, in the projective line $\PP^1$ `at infinity'.
In more detail, note we have  $\Proj(S^\bullet(H)) = \PP(H^\vee) \cong \PP(H) =
\PP^1$,
where we identify
$H^\vee = H$ via $\omega$.
Heuristically, one may view the  graded algebra morphism 
%
\begin{equation*}
A^\tau\onto A^\tau/\la z\ra\cong S^\bullet(H) \cong \BC[x,y],
\end{equation*}
as being induced by a `closed embedding' $\PP^1\into \BPt$',
of the projective line `at infinity'.
Specifically, there is a  pair of
adjoint  functors $i_*: \mathsf{coh}(\BP(H)) \to \mathsf{coh}(\BPt)$ and
$i^*: \mathsf{coh}(\BPt) \to \mathsf{coh}(\BP(H))$.
The  pushforward functor $i_*$ extends a graded $S^\bullet(H)$-module structure to a graded $A^\tau$-module structure
by setting the action of $z$ to be zero.
The  pullback functor $i^*$  takes a graded $A^\tau$-module $M$ to $M/Mz$.
The projection $A^\tau\onto S^\bullet(H)$ is clearly $\SL(H)$-equivariant, hence so are the functors $i_*$ and $i^*$.
The  functor $i_*$ is exact. The functor $i^*$ is right exact, and it has
 a sequence of left derived functors $L_pi^*$, $p > 0$. In fact, we have:
\begin{itemize}
\item for any sheaf $E$, one has $L_{>1}i^*E = 0$;
\item for a torsion free sheaf $E$, one also has $L_1i^*E = 0$;
\item for a locally free sheaf $\CE$, the sheaf $i^*\CE$ is also locally free.
\end{itemize}

We will use the following result.

\prop{morph2}(\cite[Proposition 3.4.14]{BGK})
For any $\tau\ne0$ one has

(1) If $i^*E=0$, then $E=0$.

(2) If $\phi\in\Hom(E,F)$ and $i^*\phi$ is an epimorphism, then $\phi$ is an epimorphism.

(3)  If $\phi\in\Hom(E,F)$ and both $i^*\phi$ and $L_1i^*\phi$ are isomorphisms, then $\phi$ is an isomorphism.

(4)  If $\phi\in\Hom(E,F), i^*\phi$ is a monomorphism and $L_1i^*F=0$, then $\phi$ is a monomorphism.

(5) A sheaf $E$ is locally free iff $L_{>0}i^*E=0$ and $i^*E$ is locally free.

\eprop


We deduce the following properties of Artin sheaves:

\prop{Art0} Let $F$ be an Artin sheaf and  $h_F(t)=n$.

(1) For sufficiently general $h \in H \subset \tH = H^0(\BPt,\CO(1))$,
 the map $h\colon F(-1)\to F$, of
right multiplication by $h$,  is an isomorphism.

(2) For any  locally free sheaf $\CE$ we have
\begin{equation*}
\dim\Ext^m(\CE,F)=\begin{cases}
0, &\text{if $m>0$}\\
nr(\CE), &\text{if $m=0$}.
\end{cases}
\end{equation*}


(3) The sheaf $F$ has a filtration
$0=F_0\subset F_1\subset\ldots\subset F_n=F$ such that $F_k/F_{k-1}=i_*\CO_{P_k}$
for some points $P_1,\dots,P_n \in \PP(H)$ on the line at infinity.
In particular, if $h_F(t)=1$, then $F\cong i_*\CO_P$
for some point $P\in\PP(H)$.
\eprop
\prf
(1) Since both $i^*F$ and $L_1i^*F$ are torsion sheaves on $\PP(H)$, the maps $i^*h$ and $L_1i^*h$ are isomorphisms for generic $h$.
Hence $h\colon F(-1) \to F$ is also an isomorphism for generic $h$ by \refp{morph2}(3).

(2) By~\refe{resol2} it is enough to consider the case $\CE = \CO(p)$ for some $p \in \ZZ$.
In this case,  the result is clear for $p \ll 0$. The result for arbitrary $p$ then follows from (1).


(3) The map $F\to i_*i^*F$ is an epimorphism. On the other hand, $i^*F$ is a nontrivial sheaf on $\PP(H)$,
hence there is an epimorphism $i^*F\to\CO_P$ for some $P\in\PP(H)$. The composition gives an epimorphism $F\to i_*\CO_P$.
Its kernel is an Artin sheaf on $\BPt$ of length $n-1$ and we can apply induction in $n$.
%
%
\epr


Let $F$ be an Artin sheaf and take an arbitrary $p \in \ZZ$. Consider the canonical map
\begin{equation*}
H^0(\BPt,F(p)) \otimes H \to H^0(\BPt,F(p+1))
\end{equation*}
induced by the embedding $H \subset \tH = H^0(\BPt,\CO(1))$.
Let $n = h_F$ be the length of $F$, so that both cohomology spaces above are $n$-dimensional.
A component of the $n$-th wedge power of the above map is a map
\begin{equation*}
\det (H^0(\BPt,F(p))) \otimes S^n H \to \det(H^0(\BPt,F(p+1))).
\end{equation*}
Its partial dualization gives a map
\begin{equation}\label{equation-support}
\det (H^0(\BPt,F(p))) \otimes \det(H^0(\BPt,F(p+1)))^\vee \to S^nH^\vee.
\end{equation}
We consider the projectivization of the right hand side as the space of degree $n$ divisors on $\PP(H)$
and denote by
\begin{equation*}
\supp(F)  \in \PP(S^nH^\vee) = S^n\PP(H^\vee)
\end{equation*}
the image of the map. The next lemma shows it is well-defined.

\lem{artin-support}
For any Artin sheaf $F$ of length $n$, the map~\eqref{equation-support} is injective and
$\supp(F)$ is well-defined. Furthermore, it is independent of the choice of $p$.
If $0 \to F_1 \to F_2 \to F_3 \to 0$ is an exact sequence of Artin sheaves, then
\begin{equation*}
\supp(F_2) = \supp(F_1) + \supp(F_3).
\end{equation*}
If $0 \ne h \in H$, then $h:F(-1) \to F$ is an isomorphism if and only if $h \not\in \supp(F)$.
\elem
\prf
Clearly, evaluation of the image of~\eqref{equation-support} on $h \in H$ is the determinant of the map $H^0(\BPt,F(p)) \to H^0(\BPt,F(p+1))$ induced by $h$.
We know that  the map is an isomorphism for generic $h$. Hence its determinant is nonzero. This means that the image of~\eqref{equation-support} is not
identically zero and proves the first claim of the lemma. It also proves the ``only if'' part of the last claim. Moreover, the ``if'' part also follows
for Artin sheaves of length~1. The additivity of the support under extensions is evident (the determinant of a block upper triangular matrix is the product
of the determinants of blocks). It follows that if $F_\bullet$ is a filtration on $F$ such that $F_k/F_{k-1} \cong i_*\CO_{P_k}$, then $\supp(F) = \sum P_k$.
In particular, $\supp(F)$ does not depend on the choice of the integer $p$. Finally, this observation and the additivity of the support also proves
the ``if'' part of the last claim in general.
\epr

%



One can also consider families of sheaves on a noncommutative
plane $\BPt$. More precisely, for any affine scheme $S$ one has  a notion of a coherent
sheaf on $S \times \BPt$ (see~\cite{ns}), which is the standard way to think about $S$-families of sheaves on $\BPt$.
Once on has the notion of a family, one can define moduli spaces of sheaves on $\BPt$ with appropriate stability conditions, which is the goal of this section.

Many results of this section hold in a more general setting
of an arbitrary Artin--Schelter algebra. 
In particular,
the Gieseker moduli space that we construct is the same as the moduli space studied by Nevins and Stafford ~\cite{ns}.
The special case of the algebra $A_\tau$ is, however, much simpler than the general case.
Therefore, we are going to present  the corresponding constructions in our special case in full detail
while proofs will  often be omitted since proofs of more general results can be found in ~\cite{ns}.

Uhlenbeck spaces have not been considered in {\em loc cit}.  Its definition and all the results involving
the Uhlenbeck space are new.
The relation between the Uhlenbeck and the Gieseker moduli spaces will be based on  a GIT construction of the latter.
That construction is different from the original construction of   the Gieseker space used in \cite{ns}.


\ssec{quive}{Coherent sheaves and quiver representations}
Let  $A^!_\tau$ be
the quadratic dual algebra  of $A^\tau$. The relations 
for $A^\tau$ show  that $A^!_\tau$ is isomorphic
to a 
{\sf twisted exterior algebra} $\Lambda^\bullet_\tau(\tH^\vee)$ of 
the vector space $\tH^\vee = H^\vee \oplus \BC\zeta$.
Specifically, writing
 $\{-,-\}$  for the anticommutator,
we have
\begin{equation*}
A^!_\tau \cong
\Lambda^\bullet_\tau(\tH^\vee) =
\BC\langle H^\vee \oplus \BC\zeta \rangle/\langle \{H^\vee, H^\vee\} = \{H^\vee,\zeta\} = \zeta^2 + \tau\omega = 0 \rangle.
\end{equation*}
The group $\SL(H)$ acts on~$A^!_\tau$ by algebra automorphisms.

Choosing a symplectic basis $\xi,\eta$ of $H^\vee$ we can rewrite the
above algebra as follows:
\begin{equation*}
A^!_\tau = \BC\langle \xi,\eta,\zeta \rangle/\langle \xi^2=\eta^2=\eta\xi+\xi\eta=\zeta\xi + \xi\zeta = \eta\zeta + \zeta\eta = \zeta^2 + \tau(\xi\eta - \eta\xi) = 0 \rangle.
\end{equation*}
Here, the grading on $A^!_\tau$ corresponds to the grading $\deg\xi=\deg\eta=\deg\zeta=1$.
%
Let $\bpt$ be the following quiver:
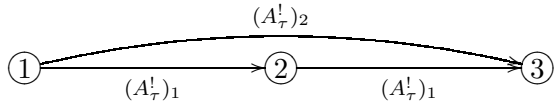
\begin{equation*}
\xymatrix{
*+[o][F]{1} \ar[rrr]_{(A_\tau^!)_1}\ar@/^1pc/[rrrrrr]^{(A_\tau^!)_2} &&&
*+[o][F]{2} \ar[rrr]_{(A_\tau^!)_1} &&&
*+[o][F]{3} }
\end{equation*}
with the spaces of arrows given by the components $(A_\tau^!)_1$ and $(A_\tau^!)_2$ of the dual algebra
and the composition of arrows given by the multiplication $(A^!_\tau)_1 \otimes (A^!_\tau)_1 \to (A^!_\tau)_2$ in $A^!_\tau$.
The $\SL(H)$ action on $A^!_\tau$ induces an action on the quiver $\bpt$,
on the category of its representations $\Rep(\bpt)$, and on its derived category $\BD(\Rep(\bpt))$.

\prop{derequ}
The functors between the bounded derived categories
\begin{align*}
&\BD(\mathsf{coh}(\BPt)) \to \BD(\Rep(\bpt)),
&&
E \mapsto (\Ext^\bullet(Q_2(-1),E), \Ext^\bullet(Q_1,E), \Ext^\bullet(Q_0(1),E)),\\
&\BD(\Rep(\bpt)) \to \BD(\mathsf{coh}(\BPt)),
&&
R_\bullet \mapsto \{ R_1\otimes\CO(-1) \to R_2\otimes\CO \to R_3\otimes\CO(1)\}
\end{align*}
are mutually inverse $\SL(H)$-equivariant equivalences.
\eprop
\prf
Follows from the fact that $(\CO(-1),\CO,\CO(1))$ is a strong exceptional collection in $\BD(\mathsf{coh}(\BPt))$,
and $(Q_2(-1),Q_1,Q_0(1))$ is its dual collection. The quiver $\bpt$ is in fact
the quiver of morphisms of the latter sequence.
\epr

We consider the restrictions of these functors to the abelian categories.
Given a representation $R_\bullet = (R_1,R_2,R_3)$ of $\bpt$, one constructs a complex of sheaves
\begin{equation}\label{comrep}
\SC(R_\bullet) := \{R_1\otimes\CO(-1) \to R_2\otimes \CO \to R_3\otimes\CO(1)\}.
\end{equation}
Denote by~$\CH^i(R_\bullet)$, $i =1,2,3$, its cohomology sheaves. Recall that a three-term complex
is a {\sf monad} if its cohomology at the first and  last terms vanish.

Similarly, given a sheaf $E$ we consider a representation of $\bpt$
\begin{equation}\label{repcom}
V_\bullet(E) = (\Ext^1(Q_2(-1),E), \Ext^1(Q_1,E), \Ext^1(Q_0(1),E)).
\end{equation}
This is equivalent to applying the functor of~\refp{derequ} and then taking
the first cohomology in the derived category of quiver representations.

\lem{monad-equivalence}
Let $\SC(R_\bu)$ be a monad and $E = \CH^2(\SC(R_\bu))$  its middle cohomology sheaf. Then $V_\bu(E) \cong R_\bu$
and $\Ext^i(Q_p(1-p),E) = 0$ for $i \ne 1$ and $p = 0,1,2$.

Conversely, let $E$ be a coherent sheaf on $\BPt$ such that $\Ext^i(Q_p(1-p),E) = 0$ for $i \ne 1$ and $p = 0,1,2$.
Then $\SC(V_\bu(E))$ is a monad and $\CH^2(V_\bu(E)) = E$.
\elem
\prf
If $\SC(R_\bu)$ is a monad, then $E = \CH^2(\SC(R_\bu))$ is isomorphic
to the complex $\SC(R_\bu)$ in the derived category $\BD(\mathsf{coh}(\BPt))$,
hence the complex can be used to compute $\Ext^i(Q_p(1-p),E)$.
The computation gives the required result.

Conversely, it follows from the assumptions of the lemma that the representation  $V_\bu(E)$ corresponds, 
under the equivalence  of ~\refp{derequ}, to the image of $E$ in $\BD(\Rep(\bpt))$.
Hence $\SC(V_\bu(E)) \cong E$ in $\BD(\mathsf{coh}(\BPt))$.
This means that the complex is a monad and its middle cohomology is $E$.
\epr




\ssec{stabilo}{Stability of sheaves and quiver representations}

The notions of Gieseker and Mumford (semi)stability of coherent sheaves are  standard in the commutative context.
We refer to~\cite{HL} for more details and for proofs of standard
facts. These notions 
have generalizations
for sheaves on $\BPt$.

Given a sheaf $E$ on $\BPt$ with $r(E)>0$,
we define its Mumford and Gieseker slopes as
\begin{align*}
\mm(E) & =\frac{\deg(E)}{r(E)}\in\BQ,\\
\mg(E) & =\frac{h_E(t)}{r(E)}=\frac{(t+1)(t+2)}{2}+\mm(E)\frac{2t+3}{2}+\frac{\deg(E)^2-2c_2(E)}{2r(E)}\in\BQ[t].
\end{align*}
%
Let $p(t)$ and $q(t)$ be polynomials. We say that $p<q$ (resp. $p\le q$) if for all $t\gg 0$, we have $p(t)<q(t)$ (resp. $p(t)\le q(t)$).

\defe{stability}
A sheaf $E$ is called {\itshape Gieseker stable} (resp.\ {\itshape Gieseker semistable})
if $E$ is torsion free and for any subsheaf $0 \subsetneq F \subsetneq E$, we have $\mg(F)<\mg(E)$ (resp.\ $\mg(F) \le \mg(E)$).
A pair of sheaves $E$ and $F$ are called Gieseker $S$-equivalent if both of them are Gieseker semistable and 
they have isomorphic composition factors in the category of Gieseker semistable sheaves.

A sheaf $E$ is called {\itshape Mumford stable} (resp.\ {\itshape Mumford semistable})
if any torsion subsheaf of $E$ is Artin and for any $F \subset E$ such that $0<r(F)<r(E)$, we have $\mm(F)<\mm(E)$ (resp.\ $\mm(F)\le\mm(E)$).
A pair of sheaves $E$ and $F$ are called Mumford $S$-equivalent if both of them are Mumford semistable and they have isomorphic composition factors in the category of Mumford semistable sheaves.
\edefe

Both Gieseker and Mumford stabilities of sheaves on $\BPt$ behave analogously to those on the commutative projective plane $\BP^2$.
For example, by~\cite{Ru} each sheaf $F$ has a {\sf Harder--Narasimhan filtration}, i.e., a filtration 
\begin{equation*}
0 = F_n \subset F_{n-1} \subset \dots \subset F_1 \subset F_0 = F
\end{equation*}
such that $F_{i-1}/F_i$ are stable and $\mu(F_{n-1}/F_n) > \mu(F_{n-2}/F_{n-1}) > \dots > \mu(F_0/F_1)$.

\medskip


To check Gieseker stability (semistability) it is sufficient to consider only subsheaves $F \subset E$ such that
$E/F$ is torsion free (in particular, $r(F) < r(E)$). So, the following is clear.

\lem{Q}
Any torsion free sheaf of rank $1$ is Gieseker stable and Mumford stable.
\elem

Note  that $\mm(E)>\mm(F)$ implies $\mg(E)>\mg(F)$, resp. $\mg(E)\ge\mg(F)$
implies $\mm(E)\ge\mm(F)$. It follows that Gieseker semistability implies
Mumford semistability, while Mumford stability for torsion free sheaves implies Gieseker stability.
Furthermore, in the case where the rank and degree of a torsion free sheaf are coprime
semistablity implies stability.

The following result is standard
\lem{hom}

(1) Let $E$, $F$ be Mumford semistable sheaves such that $F$ torsion free and $\mu_M(E)>\mu_M(F)$. Then $\Hom(E,F)=0$.

(2) Let $E,F$ be Mumford stable sheaves such that $E$ is locally free and $\mu_M(E) \ge \mu_M(F)$.
Then any nontrivial homomorphism $E\to F$ is an isomorphism.
%
\elem


The notion of stability for a representation of a quiver depends on a choice of a polarization, see~\cite{K}.
A {\sf polarization} in case of the quiver $\bpt$ amounts to a triple $\tet=(\tet_1,\tet_2,\tet_3)$ of real numbers.
The {\sf $\tet$-slope} of a representation $R_\bullet = (R_1,R_2,R_3)$ of $\bpt$ is defined as
\begin{equation*}
\mt(R_\bullet) = \langle \tet, \dim R_\bu \rangle := \tet_1\dim R_1+\tet_2\dim R_2+\tet_3\dim R_3.
\end{equation*}

\defe{stability-reps}
A representation $R_\bullet$ is called $\tet$-stable
(resp. $\tet$-semistable) if $\mt(R_\bullet)=0$ and, for any
subrepresentation $R'_\bullet\subset R_\bullet$ such that $0\ne
R'_\bullet\ne R_\bullet$, we have $\mt(R'_\bullet)>0$
(resp. $\mt(R'_\tet)\ge 0$). Representations $R_\bullet$ and
$R'_\bullet$ are called $S$-equivalent with respect to $\tet$ if both of
them are $\tet$-semistable and have isomorphic composition factors in the
category of $\tet$-semistable representations.
\edefe


Let $\tet,\tet'$ be a pair of polarizations. It is well-known (e.g.,~\cite{DH})
that, for all sufficiently small and positive $\eps\in \mathbb{R}$,
 stability, semistability and $S$-equivalence with respect to
 $\tet+\eps\tet'$ does not depend on $\eps$.

\defe{stability-reps-2}
A representation $R_\bullet$ is $\ttt$-stable (resp. $\ttt$-semistable) if $R_\bullet$ is
$(\tet+\eps\tet')$-stable (resp. semistable) for sufficiently small positive
$\eps$.
\edefe

There is an analogue of~\refl{hom} for representations of the quiver $\bpt$.



\ssec{sheaves-reps}{From sheaves to quiver representations}
The following result is essentially a combination of Lemma~6.4 and Theorem~5.6 from~\cite{ns}.
The only new statement is the exactness claim. We provide a proof for
the
reader's convenience.
\th{mon}
Let $-r \le d < r$.
Then, the assignment $E\mto V_\bu(E)$ gives an exact functor  from the
category of Mumford semistable torsion free sheaves $E$
on $\BPt$ such that  $r(E)=r$ and $\deg(E)=d$
to the category of representations
of the quiver $\bpt$.
For such a  sheaf $E$,
the representation $V_\bu(E)$ gives  a monad
\begin{equation}\label{monad-of-e}
V_1(E) \otimes \CO(-1) \to V_2(E) \otimes \CO \to V_3(E) \otimes \CO(1)
\end{equation}
such that its cohomology at the middle term is isomorphic to
$E$. Furthermore,  we have   $\dim V_\bullet(E) =
(n-d(d-1)/2,2n-d^2+r,n-d(d+1)/2)$,
where $n=c_2(E)$; in particular, $c_2(E) \ge d(d+1)/2$.
\eth
\prf
%
First we note that all $Q_p$ are Mumford stable of slopes equal to $0$, $3/2$ and $3$ respectively.
Indeed, for $p = 0$ and $p = 2$ this follows from~\refl{Q} and the definitions of $Q_p$.
So let $p = 1$. The sheaf $Q_1$ is locally free because by~\eqref{anotherq1} it is the kernel
of a morphism of locally free sheaves, and moreover $r(Q_1) = 2$, $\deg(Q_1) = 3$.
So, it is enough to check that if $F \subset Q_1$ is a subsheaf of rank 1 with $Q_1/F$ torsion free, then $\deg(F) \le 1$.
Assume $\deg(F) \ge 2$. Since $Q_1/F$ is torsion free, $F$ is locally free.
Since $Q_1$ is a subsheaf in $\calO(2)^{\oplus3}$ there is a nontrivial homomorphism from $F$ to $\calO(2)$.
On the other hand, both $F$ and $\calO(2)$ are Mumford stable by (1), $F$ is locally free, and
$\mm(F)\ge2=\mm(\calO(2))$. Hence $F\simeq\CO(2)$ by \refl{hom}(2). But applying
the functor $\Hom(\calO(2),-)$ to~\eqref{anotherq1} we see that $\Hom(\calO(2),Q_1)=0$.

The proved stability implies that
\begin{equation*}
\Hom(Q_0(1),E) = \Hom(Q_1,E) = \Hom(Q_2(-1),E) = 0.
\end{equation*}
Indeed, the slopes of the first arguments are $1$, $3/2$, and $2$ respectively, while the slope
of the second argument is $d/r < 1$, so~\refl{hom}(1) applies. Analogously,
\begin{equation*}
\Hom(E(3),Q_0(1)) = \Hom(E(3),Q_1) = \Hom(E(3),Q_2(-1)) = 0
\end{equation*}
since the slope of the first argument is $d/r+3 > 2$. By Serre duality we then have
\begin{equation*}
\Ext^2(Q_0(1),E) = \Ext^2(Q_1,E) = \Ext^2(Q_2(-1),E) = 0.
\end{equation*}
Therefore, \refl{monad-equivalence} applies to $E$ and shows that~\eqref{monad-of-e} is a monad
and $E$ is its cohomology.
%
The dimensions of the spaces $V_p(E)$ are computed directly by using the formula~\refe{david} for the Hilbert polynomial of a sheaf.
The exactness of the functor $V_\bu$ is clear from its definition and vanishing of $\Hom$ and $\Ext^2$ spaces.
\epr

\prop{artin-mon}
The functor $F\mto {\mathcal
  C}(R_\bullet(F))$
yields, for an Artin sheaf $F$,  a canonical exact sequence
\begin{equation*}
0 \to W_1(F) \otimes \CO(-1) \to W_2(F) \otimes \CO \to W_3(F) \otimes \CO(1) \to F \to 0.
\end{equation*}
The resulting functor $W_\bu$ from the category of Artin sheaves on $\BPt$ to the category
of representations of the quiver $\bpt$ is exact
and we have  $\dim W_\bullet(F) = (n,2n,n)$, where $n$ is the length of $F$.
\eprop
\prf
The proof is analogous to the proof of~\refl{monad-equivalence}. We apply the equivalence of~\refp{derequ} to the sheaf~$F$.
By~\refp{Art0}(2), applying the functor of~\refp{derequ} to $F$ yields the representation
\begin{equation*}
W_\bullet(F) = (\Hom(Q_2(-1),F),\Hom(Q_1,F),\Hom(Q_0(1),F)).
\end{equation*}
Its dimension vector equals $(n,2n,n)$. Since the functor is an equivalence, it follows that
the complex $\SC(W_\bu(F))$ is left exact and $\CH^3(\SC(W_\bullet(F))) \cong F$, which amounts to the above exact sequence.
Exactness of the functor $W_\bu$ follows from the vanishing of $\Ext^1(Q_p(1-p),F)$ by~\refp{Art0}(2).
\epr

If $0 \ne h \in H$, $P \in \PP(H)$ is the corresponding point and $F = \CO_P$, then
\begin{equation}\label{monad-point}
W_\bullet(\CO_P) = \{ \C \xrightarrow{\ \left(\begin{smallmatrix} h \\ \zeta \end{smallmatrix}\right)\ } \C^2 \xrightarrow{\ (-\zeta,h)\ } \C \}.
\end{equation}


\ssec{relat stab}{From sheaf stability to quiver stability}

In this section we show that Gieseker and Mumford semistability correspond to the semistablity of quiver representations.

From now on we fix a triple $(r,d,n)$ such that
\begin{equation}\label{assumptions-r-d-n}
0 \le d < r
\qquad\text{and}\qquad
n \ge d(d+1)/2.
\end{equation} Put
\begin{equation}\label{alpha}
\alp(r,d,n)=(n-d(d-1)/2, 2n-d^2+r, n-d(d+1)/2).
\end{equation}
Let $E$ be a Mumford semistable sheaf such that $r(E)=r,
\deg(E)=d$ and $c_2(E)=n$. Then $\dim V_\bullet(E) = \alp(r,d,n)$ by \reft{mon}.

We choose the following pair of polarizations:
\begin{equation}\label{thetas}
\begin{split}
&\tet^0=(-r-d,d,r-d),\\
&\tet^1=(2n-d^2+r,d^2-2n,2n-d^2+r).
\end{split}
\end{equation}
Note that $\tet^0$ does not depend on $n$. Note also that
\begin{equation*}
\langle \tet^0,\alp(r,d,n) \rangle = \langle \tet^1,\alp(r,d,n) \rangle = 0.
\end{equation*}
In what follows we frequently consider $\tet^0$-stability and $\tnt$-stability of representations.
In fact the notion of $\tnt$-(semi)stablity of a quiver representation is equivalent to the notion
of (semi)stability of a Kronecker complex considered in \cite{ns}.

\lem{kronecker-stability}
Let $V_\bu$ be an $\alp(r,d,n)$-dimensional representation of the quiver $\bpt$ and let $\SC(V_\bu)$
be the associated complex~\eqref{comrep}. Then $V_\bu$ is $\tnt$-(semi)stable if and only if $\SC(V_\bu)$
is (semi)stable in the sense of~\cite[Def.~6.8]{ns}.
\elem
\prf
Just note that a subcomplex in $\SC(V_\bu)$ always corresponds to a subrepresentation $U_\bu \subset V_\bu$,
and the expression (1) from \cite[Defenition~6.8]{ns} for $\SC(U_\bu)$ equals $\mu_{\tet^0}(U_\bu)$,
while the expression (2) equals $\mu_{\tet^1}(U_\bu)$,
\epr

The following crucial observation that relates the notion of stability for sheaves and for
quiver representations, respectively,
is due to Sergey Kuleshov. Parts (2) and (3) (as well as a part of~\refl{monadic-sheaf} below)
are also proved in Proposition~6.20 of~\cite{ns}.

\lem{eq}
Let $E$ be a torsion free sheaf with $r(E)=r, \deg(E)=d, c_2(E)=n$ and let
$V_\bullet(E)$ be the corresponding representation of the quiver $\bpt$. Then

(1) if $E$ is Mumford semistable, then $\VE$ is $\tet^0$-semistable;

(2) if $E$ is Gieseker semistable, then $\VE$ is $\tnt$-semistable;

(3) if $E$ is Gieseker stable, then $\VE$ is $\tnt$-stable.
\elem
%

\prf
(1) Assume that $E$ is Mumford semistable. Let $U_\bullet$ be a subrepresentation
of $\VE$, let $W_\bullet = V_\bullet(E) / U_\bullet$  be the quotient representation,
and put $u_i=\dim U_i$. We have a short exact sequence of complexes
\begin{equation}\label{uvw}
\vcenter{\xymatrix{
0 \ar[r] &
U_1 \otimes \CO(-1) \ar[r] \ar[d] &
U_2 \otimes \CO \ar[r] \ar[d] &
U_3 \otimes \CO(1) \ar[r] \ar[d] &
0
\\
0 \ar[r] &
V_1(E) \otimes \CO(-1) \ar[r] \ar[d] &
V_2(E) \otimes \CO \ar[r] \ar[d] &
V_3(E) \otimes \CO(1) \ar[r] \ar[d] &
0
\\
0 \ar[r] &
W_1 \otimes \CO(-1) \ar[r] &
W_2 \otimes \CO \ar[r] &
W_3 \otimes \CO(1) \ar[r] &
0.
}}
\end{equation}
We  view the above diagram  as an exact triple (with respect to the vertical maps) of 3-term complexes
and apply the Snake Lemma. We obtain a long exact sequence of cohomology sheaves:
\begin{equation*}
0 \to
\CH^1(U_\bullet) \to 0 \to \CH^1(W_\bullet) \to
\CH^2(U_\bullet) \to E \to \CH^2(W_\bullet) \to
\CH^3(U_\bullet) \to 0 \to \CH^3(W_\bullet) \to 0
\end{equation*}
of these complexes. In particular, we have $\CH^1(U_\bullet) = \CH^3(W_\bullet) = 0$.
%
Put
\begin{equation*}
r_i^U = r(\CH^i(U_\bullet)), \quad
d_i^U = \deg(\CH^i(U_\bullet)), \quad
r_i^W = r(\CH^i(W_\bullet)), \quad
d_i^W = \deg(\CH^i(W_\bullet)).
\end{equation*}
Let $I$ be image of the morphism $\CH^2(U_\bullet) \to E$ in the above sequence
and let
\begin{equation*}
r_I=r(I), \quad d_I=\deg(I).
\end{equation*}
Then using additivity of rank and degree we can rewrite the slope of $U_\bullet$ as
\begin{multline*}
\mu_{\tet^0}(U_\bullet)
 = r(u_3-u_1)+d(u_2-u_1-u_3)
 =r(d_3^U - d_2^U)+d(r_2^U-r_3^U) \\
 =r(d_3^U - d_1^W - d_I)+d(r_1^W+r_I-r_3^U)
=(rd_3^U - r_3^U d)+(r_1^W d-r d_1^W)+(r_Id-rd_I).
\end{multline*}
Now we will show that all three summands on the right-hand-side are nonnegative.

First, note that $\CH^3(U_\bullet)$ is a quotient of $U_3\otimes\CO(1)$.
The latter sheaf is semistable by \refl{Q}(1) and
we have $\mm(U_3\otimes\calO(1))=1$. Therefore, $d_3^U \ge r_3^U$ and hence
\begin{equation*}
rd_3^U-r_3^Ud \ge r_3^U(r-d)\ge0.
\end{equation*}
Note that the inequality $rd_3^U-r_3^Ud \ge0$ is strict unless $r_3^U=d_3^U=0$,
that is unless $\CH^3(U_\bullet)$ is an Artin sheaf.

Further, note that $\CH^1(W_\bullet)$ is a subsheaf of $W_1\otimes\CO(-1).$
The latter sheaf is semistable and one has $\mm(W_1\otimes\CO(-1))=-1$. Therefore,
$d_1^W \le -r_1^W \le 0$ and hence
\begin{equation*}
r_1^W d-r d_1^W \ge dr_1^W \ge 0.
\end{equation*}
Note that this inequality is strict unless $r_1^W = d_1^W = 0$, that is unless
$\CH^1(W_\bullet)$ is an Artin sheaf. But since it is a subsheaf in $W_1\otimes\CO(-1)$
it is torsion free, hence this is equivalent to $\CH^1(W_\bullet) = 0$.


Finally, $I$ is a subsheaf of the Mumford semistable torsion free sheaf $E$.
Hence either $I=0$, or else $r_I>0$ and $\mm(I)\le\mm(E)$. In both cases we have
$$r_Id-d_Ir\ge0.$$
Note that this inequality is strict unless $I=0$ or $\mm(I)=\mm(E).$

Combining all these inequalities we see that any subrepresentation in $\VE$
has a non-negative $\tet^0$-slope. Thus $\VE$  is $\tet^0$-semistable.

(2) Assume that the sheaf $E$ is Gieseker semistable but $V_\bullet(E)$ is not $(\tet^0,\tet^1)$-semistable.
Let $U_\bullet \subset V_\bullet(E)$ be a destabilizing subrepresentation. The sheaf
$E$ is automatically Mumford semistable. Hence $V_\bullet(E)$ is $\tet^0$-semistable by (1).
We conclude that $\mu_{\tet^0}(U_\bullet) = 0$. As  in the above argument, 
we deduce that the sheaf $\CH^3(U_\bullet)$ is an Artin sheaf and $\CH^1(W_\bullet) = 0$.
Thus $\CH^2(U_\bullet) = I$ is a subsheaf in $E$.

Let $c = \cch_2(E) = d^2/2 - n$ and $c_i^U = \cch_2(\CH^i(U_\bullet))$.
Then
\begin{multline*}
\mu_{\tet^1}(U_\bullet)=(r-2c)u_1+2cu_2+(r-2c)u_3=r(u_1+u_3)+2c(u_2-u_1-u_3) \\
=2r(c_3^U-c_2^U) + 2c(r_2^U-r_3^U)
=2(rc_3^U-r_3^U c)+2(r_2^Uc-r c_2^U).
\end{multline*}
Since $\CH^3(U_\bullet)$ is an Artin sheaf we have $r_3^U = 0$ and $c_3^U \ge 0$,
hence
\begin{equation*}
rc_3^U-r_3^U c = rc_3^U \ge 0.
\end{equation*}
On the other hand, $\CH^2(U_\bullet) = I$ and hence
either $\CH^2(U_\bullet) = 0$, or $\mm(\CH^2(U_\bullet)) = \mm(E)$.
In both cases
\begin{equation*}
r_2^Uc-r c_2^U = rr_2^U(\mg(E) - \mg(\CH^2(U_\bullet))) \ge 0
\end{equation*}
because $E$ is Gieseker semistable.
%
%
%
We deduce that $\mu_{\tet^1}(U_\bullet) \ge0$. This contradicts the assumption that $U_\bullet$
destabilizes $V_\bullet(E)$.

(3) In the notation of (2) assume that $\mu_{\tet^0}(U_\bu) = \mu_{\tet^1}(U_\bu) = 0$.
Then $c_3^U=0$, hence $\CH^3(U_\bullet)=0$. Moreover,
either $\CH^2(U_\bullet)=0$, or $\mg(\CH^2(U_\bullet))=\mg(E)$, and so $\CH^2(U_\bullet)=E$
since $E$ is Gieseker stable. In the first case the first line of~\eqref{uvw} is exact,
hence $U_\bullet=0$ by~\refp{derequ}. In the second case the first line of~\eqref{uvw} is a resolution of $E$,
hence $U_\bullet=\VE$.
\epr

\lem{cfart}
Let $F$ be an Artin sheaf. Then $W_\bu(F)$ is $\tet^0$-semistable. If the length of $F$ equals $1$, then $W_\bu(F)$ is $\tet^0$-stable.
\elem
\prf
Since by~\refp{Art0}(3) any Artin sheaf is an extension of the structure sheaves of points and the functor $W_\bullet$
is exact, it is enough to verify the $\tet^0$-stability of $W_\bullet(\CO_P)$. The latter is clear from the explicit
form~\eqref{monad-point} of the monad --- it is easy to see that the only nontrivial subrepresentations of $W_\bu(\CO_P)$
have dimension $(0,0,1)$, $(0,1,1)$ and $(0,2,1)$, and their $\tet^0$-slope is clearly positive with our assumptions on $d$ and $r$.
%
\epr

\ssec{quiver-sheaf}{From quiver stability to sheaf stability}

In this section we show that stable representations of the quiver, in their turn,
give rise to stable sheaves.

\defe{red}
A representation $V_\bullet$ is called:
\begin{itemize}
\item
{\sf Artin}, if $\CH^1(V_\bu) = \CH^2(V_\bu) = 0$ and $\CH^3(V_\bu)$ is an Artin sheaf;
\item
{\sf monadic}, if $\CH^1(V_\bullet) = \CH^3(V_\bullet) = 0$;
\item
{\sf supermonadic}, if both $V_\bu$ and $V_\bu^\vee$ are monadic.
\end{itemize}
\edefe

\lem{smlf}
A monadic representation $V_\bullet$ is supermonadic iff $\CH^2(V_\bullet)$ is locally free.
\elem
\prf
Since $V_\bullet$ is monadic, the complex $\SC(V_\bullet)$ is isomorphic to $\CH^2(V_\bullet)$ in the derived category $\BD(\mathsf{coh}(\BPt))$.
Therefore the complex $\SC(V_\bullet^\vee) = \SC(V_\bullet)^*$ is isomorphic to the derived dual of $\CH^2(V_\bu)$. In other words,
$\CH^i(\SC(V^\vee_\bu)) \cong \underline{\Ext}^{i-2}(\CH^2(V_\bullet),\CO)$. So, $V_\bullet$ is supermonadic iff
$\underline{\Ext}^i(\CH^2(V_\bullet),\CO) = 0$ for $i \ne 0$, i.e., iff $\CH^2(V_\bullet)$ is locally free.
\epr

\lem{monadic-sheaf}
Let $V_\bu$ be a monadic representation of $\bpt$ and let $E = \CH^2(V_\bu)$.
If $V_\bu$ is $\tet^0$-semistable, then $E$ is Mumford semistable.
If $V_\bu$ is $\tnt$-semistable, then $E$ is Gieseker semistable.
In each case we have $V_\bu = V_\bu(E)$.
\elem
\prf
Assume $E$ is not Mumford semistable and consider its Harder--Narasimhan filtration.
Breaking it up at slope $\mu_M(E)$, we can represent $E$ as an extension
\begin{equation*}
0 \to E' \to E \to E'' \to 0,
\end{equation*}
such that the slopes of all quotients in the Harder--Narasimhan filtration of $E'$ (resp.\ $E''$) are greater than (resp.\ less than or equal to) $\mu_M(E)$.
Let $(r',d',n')$
be the rank, the degree and the second Chern class of $E'$.
Note that both $E'$ and $E''$ are the cohomology sheaves of the monads $V_\bu(E')$ and $V_\bu(E'')$ respectively.
Indeed, for $E'$ the argument of~\reft{mon} shows that $\Ext^2(Q_p(1-p),E') = 0$.
On the other hand, since $E'$ is a subsheaf of $E$, we have $\Hom(Q_p(1-p),E') \subset \Hom(Q_p(1-p),E) = 0$.
Similarly, for $E''$ the argument of~\reft{mon} gives the vanishing of $\Hom$'s, while the surjectivity
of the map from $E$ gives the vanishing of $\Ext^2$. It follows that we have an exact sequence of monads
\begin{equation*}
0 \to V_\bu(E') \to V_\bu \to V_\bu(E'') \to 0.
\end{equation*}
Finally, note that
\begin{equation*}
\langle \tet^0, \alpha(r',d',n') \rangle = r'd - rd' = rr'\left(\frac{d}{r} - \frac{d'}{r'}\right) = rr'(\mu_M(E) - \mu_M(E')) < 0.
\end{equation*}
Hence the subrepresentation $V_\bu(E') \subset V_\bu$ violates
the $\tet^0$-semistability of $V_\bu$. This proves the first part.

If $E$ is Mumford semistable but not Gieseker semistable, we again take $E'$ to be the part of the Harder--Narasimhan
filtration of $E$ with the slopes greater than $\mu_G(E)$. Then $\langle \tet^0, \alpha(r',d',n') \rangle = 0$ but
\begin{equation*}
\langle \tet^1, \alpha(r',d',n') \rangle = r'(d^2-2n) - r(d'^2 - 2n') = 2rr'\left(\frac{d^2-2n}{2r} - \frac{d'^2-2n'}{2r'}\right) < 0.
\end{equation*}
Hence the subrepresentation $V_\bu(E') \subset V_\bu$ violates the $\tnt$-semistability of $V_\bu$. This proves the second part.

Finally, we have $V_\bu = V_\bu(E)$ by~\refl{monad-equivalence}.
\epr

\prop{tet0-ss}
Let $V_\bu$ be a $\tet^0$-semistable representation of $\bpt$ of dimension $\alpha(r,d,n)$.
Then $V_\bu$ is S-equivalent to a direct sum $U_\bu \oplus W_\bu$, where $U_\bu$ is supermonadic of dimension $\alpha(r,d,n-k)$
and $W_\bu$ is Artin of dimension $(k,2k,k)$ for some $0 \le k \le n$.
\eprop

\prf
Assume that $\CH^3(V_\bu) \ne 0$. Then $i^*(\CH^3(V_\bullet)) \ne 0$ by~\refp{morph2}(1), hence
there is a surjective morphism $\CH^3(V_\bullet) \to i_*\CO_P$ for some point $P \in \BP(H)$.
Since $\CH^3(V_\bullet)$ is the top cohomology of the complex $\SC(V_\bullet)$, there is
a canonical morphism $\SC(V_\bullet) \to \CH^3(V_\bullet)$.
Composing these morphisms we get a nontrivial morphism $\SC(V_\bullet) \to \CO_P$.
By \refp{derequ} this corresponds to a nontrivial morphism $V_\bullet \to W_\bullet(\CO_P)$.
Since $V_\bu$ is $\tet^0$-semistable and $W_\bu(\CO_P)$ is $\tet^0$-stable, the morphism is surjective.
Taking $V'_\bu$ to be the kernel of the morphism, we see that $V'_\bu$ is $\tet^0$-semistable and $V_\bu$ is S-equivalent to
$V'_\bu \oplus W_\bu(\CO_P)$. The dimension of $V'_\bu$ is strictly less than that of $V_\bu$,
so iterating the construction we reduce to the case when $\CH^3(V_\bu) = 0$.


Assume now that $\CH^3(V_\bu) = 0$, but $\CH^3(V^\vee_\bu) \ne 0$. Then applying the same argument to $V^\vee_\bu$
we obtain an injection $W_\bu(\CO_P) \to V_\bu$. Taking $V'_\bu$ to be the cokernel of this morphism,
we see that $V'_\bu$ is $\tet^0$-semistable and $V_\bu$ is S-equivalent to $V'_\bu \oplus W_\bu(\CO_P)$.
Iterating the construction we reduce to the case when $\CH^3(V_\bu) = \CH^3(V^\vee_\bu) = 0$.

Finally, assume that $\CH^3(V_\bu) = 0$ and $\CH^3(V^\vee_\bu) = 0$. Let us show that $V_\bu$ is supermonadic.
Let $\CF := \CH^1(V_\bu) $ and assume $\CF\ne 0$. Then $\CF$ is a locally free sheaf and we have a left exact sequence
$0 \to \CF \to V_1 \otimes \CO(-1) \to V_2 \otimes \CO$. Applying dualization we get a complex
\begin{equation*}
V_2^\vee \otimes \CO \to V_3^\vee \otimes \CO(1) \to \CF^*.
\end{equation*}
The second arrow here is nontrivial since it induces the embedding $\CF \to V_1\otimes\CO(-1)$
after the second dualization.
It follows  that $\CH^3(V^\vee_\bu) \ne 0$, contradicting the assumption.
We conclude that $\CH^1(V_\bu) = 0$.
An similar argument with $V_\bu$ replaced by $V^\vee_\bu$ shows that $\CH^1(V^\vee_\bu) = 0$.
Hence $V_\bu$ is  supermonadic.

At each step of the above procedure the dimension of the representation
decreases by $(1,2,1)$. Therefore, at the end of the procedure, 
the dimension of the supermonadic part is equal to
\begin{equation*}
\alpha(r,d,n) - (k,2k,k) = \alpha(r,d,n-k).
\end{equation*}
Since all the components of a dimension vector are nonnegative, we have $n-k \ge d(d+1)/2$; in particular,
we get  $k \le n$.
\epr
Note that we have
\[(n,2n,n) = \alpha(0,0,n).\]
\cor{tet-ss-artin}
Let $W_\bu$ be a $\tet^0$-semistable representation of $\bpt$ of dimension $(n,2n,n)$. Then $W_\bu$ is an Artin representation.
\ecor
\prf
Applying~\refp{tet0-ss} we see that $W_\bu$ is S-equivalent to a sum of
an Artin representation and a supermonadic representation $U_\bu$ of dimension $\alpha(0,0,n-k)$ for some $k$.
By definition and~\refl{smlf} the corresponding complex $\SC(U_\bu)$ is a monad and its middle cohomology
is a locally free sheaf of rank $0$. Thus the cohomology is zero and the complex $\SC(U_\bu)$ is acyclic.
By~\refp{derequ} this means that $U_\bu = 0$, so we have no supermonadic part. It follows that $W_\bu$
is S-equivalent to an Artin representation. It follows immediately that $\CH^1(W_\bu) = \CH^2(W_\bu) = 0$,
and $\CH^3(W_\bu)$ is an iterated extension of Artin sheaves. Hence it is an Artin sheaf itself,
and so $W_\bu$ is also an Artin representation.
\epr

The first part of the following result can be found in Lemma~6.14 of \cite{ns}.

\prop{tnt-ss}
Let $V_\bu$ be a $\tnt$-semistable representation of $\bpt$ of dimension $\alpha(r,d,n)$ with $0 \le d < r$.
Then $V_\bu$ fits into a short exact sequence
\begin{equation*}
0 \to W_\bu \to V_\bu \to U_\bu \to 0,
\end{equation*}
where $U_\bu$ is supermonadic and $W_\bu$ is Artin.
Moreover, $V_\bu = V_\bu(E)$, where $E$ is a Gieseker semistable sheaf of rank $r$, degree $d$ and $c_2 = n$,
$U_\bu = V_\bu(E^{**})$, and $W_\bu = W_\bu(E^{**}/E)$.
\eprop

\prf
The argument of~\refp{tet0-ss} proves that there is a filtration on $V_\bu$ in which there are several factors
which are Artin representations of dimension $(1,2,1)$ and one supermonadic factor of dimension $\alpha(r,d,n-k)$
for some $0 \le k \le n$. But
\begin{equation*}
\langle \tet^1, (1,2,1) \rangle = 2r > 0,
\end{equation*}
hence Artin factors can appear only before the supermonadic factor. This proves that the filtration
gives the required exact sequence.

We apply the functor $\SC$ to the exact sequence and take into account that $\CH^2(W_\bu) = \CH^3(U_\bu) = 0$.
We get the long exact sequence of cohomology sheaves
\begin{equation*}
0 \to \CH^2(V_\bu) \to \CH^2(\CU_\bu) \to \CH^3(W_\bu) \to \CH^3(V_\bu) \to 0.
\end{equation*}
Assume  $\CH^3(V_\bu) \ne 0$, then t
The argument of the proof of~\refp{tet0-ss} shows that if $\CH^3(V_\bu) \ne 0$, then 
there is a surjection $V_\bu \to W_\bu(\CO_P)$.  This contradicts $\tnt$-semistability of $V_\bu$.
Thus $V_\bu$ is monadic. Put $E = \CH^2(V_\bu)$, $\CE = \CH^2(U_\bu)$ and $F = \CH^3(W_\bu)$. Then the above sequence
takes the form:
\begin{equation*}
0 \to E \to \CE \to F \to 0.
\end{equation*}
Note that $\CE$ is locally free by~\refl{smlf}, and $F$ is Artin since $W_\bu$ is. Dualizing the sequence and taking into account
that $\underline{\Hom}(F,\CO) = \underline{\Ext}^1(F,\CO) = 0$ since $F$ is Artin,
we deduce that $E^* = \CE^*$. Therefore $E^{**} = \CE^{**} = \CE$ since $\CE$ is locally free and the map $E \to \CE = E^{**}$
is the canonical embedding. Thus $E$ is torsion free and $F \cong E^{**}/E$. Moreover, by~\refp{derequ}
it follows that $V_\bu = V_\bu(E)$ and $U_\bu = V_\bu(E^{**})$, while $W_\bu = W_\bu(E^{**}/E)$.

Finally,  $E$ is Gieseker semistable by~\refl{monadic-sheaf}.
\epr

\cor{coprime-stable}
Assume that $r$ and $d$ are coprime. Then a $\tet^0$-semistable representation is $\tnt$-semistable if and only if it has no Artin quotients,
i.e.,  it is monadic.
\ecor
\prf
Let  $r$ and $d$ be coprime and
let $V_\bu$ be a supermonadic $\tet^0$-semistable representation. We claim  that  $V_\bu$  is $\tet^0$-stable.
Indeed, we have $V_\bu = V_\bu(E)$ for a Mumford semistable sheaf $E$ by~\refl{monadic-sheaf}. Moreover, the sheaf $E$ is locally free
by~\refl{smlf}. So, if $0 \to V'_\bu \to V_\bu \to V''_\bu \to 0$ is an exact sequence of representations where both $V'_\bu$ and $V''_\bu$ nonzero and
$\mu_{\tet^0}(V'_\bu) = \mu_{\tet^0}(V''_\bu) = 0$, then by~\refl{eq} we have an exact sequence
\begin{equation*}
0 \to \CH^2(V'_\bu) \to E \to \CH^2(V''_\bu) \to \CH^3(V'_\bu) \to 0.
\end{equation*}
Moreover, $V''_\bu$ is monadic, $\mu_M(\CH^2(V''_\bu)) = \mu_M(E) = d/r$, and $\CH^3(V'_\bu)$ is Artin.
Since $r$ and $d$ are coprime, we have that either the rank and the degree of $\CH^2(V''_\bu)$ are equal to zero, or 
they are equal to $r$ and $d$, respectively.
The first case is impossible since then $F := \CH^2(V''_\bu)$ is an Artin sheaf and $V''_p = \Ext^1(Q_p(1-p),F) = 0$, so $V''_\bu = 0$.
In the second case  the rank of $\CH^2(V'_\bu)$ equals zero and, hence, we must have  $\CH^2(V'_\bu)=0$ since
this sheaf is a subsheaf of $E$.
We conclude that $V'_\bu$ is a nonzero Artin subrepresentation of $V_\bu$ which means that $V^\vee_\bu$
has a nonzero Artin quotient representation and hence cannot be monadic.

Now let $V_\bu$ be a $\tet^0$-semistable representation. 
In order to prove that  $V_\bu$ is $\tnt$-semistable we must show that, for
any subrepresentation  $V' _\bu\subset V_\bu$ such that $\theta_0(V'_\bu) = 0$, we have $\theta_1(V'_\bu) > 0$.
To this end, we choose a composition series of  $V_\bu$  
in the category of $\tet^0$-semistable representationsn with zero  $\tet^0$-slope.
According to  ~\refp{tet0-ss}, all composition factors of the filtration but one are Artin representation
and the remaining composition factor, $U$, is a supermonadic representation
such that  $\theta^1(U) < 0$.
For each Artin composition factor $W_i$ we have $\theta^1(W_i) > 0$ and
$\theta^1(U + \sum W_i) = 0$. Therefore,  $V_\bu$ is $\tnt$-semistable iff 
there is a surjection $V_\bu\onto U$ with kernel $W$ which is an Artin subrepresentation such that
$\theta^1(W)>0$. We deduce that $V_\bu$ has no Artin quotients iff $V_\bu=W$, that is, iff $V_\bu$ is
a supermonadic representation.
\epr

\ssec{moduli-spaces}{Moduli spaces}

Let $\CM^\tet_\tau(r_1,r_2,r_3)$ denote the moduli space of $\tet$-semistable
$(r_1,r_2,r_3)$-dimensional representations of the quiver $\bpt$, as
defined by
King ~\cite{K}.
It is a coarse moduli space for families of $\tet$-semistable representations
of the quiver $\bpt$ of dimension $(r_1,r_2,r_3)$. In particular, its closed points
are in  bijection with S-equivalence classes of $\tet$-semistable representations.

For rational $\tet$ there is an explicit GIT construction of the moduli space. One starts with {\sf the representation space} of $\bpt$:
\begin{equation}\label{rt-embedded}
\Rt(r_1,r_2,r_3) \subset \Hom(\C^{r_1}\otimes (A_\tau^!)_1,\C^{r_2}) \times \Hom(\C^{r_2}\otimes (A_\tau^!)_1,\C^{r_3}),
\end{equation}
consisting of those pairs of maps $f:\C^{r_1} \otimes (A_\tau^!)_1 \to \C^{r_2}$, $g:\C^{r_2} \otimes (A_\tau^!)_1 \to \C^{r_3}$ such that the composition
$g\circ (f \otimes \id):\C^{r_1} \otimes (A_\tau^!)_1 \otimes (A_\tau^!)_1 \to \C^{r_3}$ factors through the multiplication map
$\C^{r_1} \otimes (A_\tau^!)_1 \otimes (A_\tau^!)_1 \to \C^{r_1} \otimes (A_\tau^!)_2$. Clearly, \eqref{rt-embedded} is a Zarisky closed subset
in an affine space.
The group
\begin{equation*}
\GL(r_1,r_2,r_3)=\GL(r_1)\times \GL(r_2)\times \GL(r_3)
\end{equation*}
acts naturally on  $\Rt(r_1,r_2,r_3)$.
Given a rational polarization $\tet$, in the trivial
bundle,
let $\C[\Rt(r_1,r_2,r_3)]^{\GL(r_1,r_2,r_3),p\tet}$
be the vector space of polynomial  $\GL(r_1,r_2,r_3)$-semiinvariants of
weight $p\tet$ (this space is  declaired to be zero  unless
$p\tet$ is an integral weight).
One defines an associated GIT quotient by
\begin{equation*}
\Rt(r_1,r_2,r_3) /\!/_\tet \GL(r_1,r_2,r_3) := \Proj\left(\bigoplus_{p = 0}^\infty \C[\Rt(r_1,r_2,r_3)]^{\GL(r_1,r_2,r_3),p\tet}\right).
\end{equation*}
Then,
according to ~\cite{K}, one has 
$\CM^\tet_\tau(r_1,r_2,r_3) \cong \Rt(r_1,r_2,r_3) /\!/_\tet
\GL(r_1,r_2,r_3)$.
Further,
it turns out that the space of all polarizations $\tet$ has a chamber structure
and the moduli space $\CM^\tet_\tau(r_1,r_2,r_3)$ depends only on the chamber in which $\tet$ sits.
This makes it possible  to define $\CM^\tet_\tau$ for an arbitrary (real) polarization $\tet$ by taking rational $\tet'$ in the same chamber as $\tet$
and setting $\CM^\tet_\tau(r_1,r_2,r_3) := \CM^{\tet'}_\tau(r_1,r_2,r_3)$.

Analogously one constructs a coarse moduli space $\Mtt(r_1,r_2,r_3)$ for a pair of polarizations $(\tet,\tet')$ by taking
an arbitrary polarization in the chamber containing $\tet + \eps\tet'$ for all sufficiently small and positive $\eps$.

It has been shown in \cite{K} that  the moduli space of
semistable representations  of any  quiver that has no
 oriented cycles is a projective variety. It follows,
 since
the quiver $\bpt$ has no oriented cycles,
that each of the above moduli spaces $\CM^\tet_\tau(r_1,r_2,r_3)$
is a projective variety. This variety comes equipped
with a  natural $\SL(H)$-action.
Finally, we remark that  if the dimension vector
$(r_1,r_2,r_3)$ is primitive, i.e.,  indivisible, then 
 $\CM^\tet_\tau(r_1,r_2,r_3)$ is a {\em fine}  moduli space.
\medskip

\nc{\reduced}{{\mathrm{red}}}

Below we discuss moduli spaces of several classes of representations of the quiver $\bpt$.
First, recall that by~\refc{tet-ss-artin} any $\tet^0$-semistable representation of dimension $(n,2n,n)$
is an Artin representation. So, we refer to the corresponding moduli space as to {\em the moduli space of Artin representations}
and denote it by $\aMt(n,2n,n)$. Thus we have
%
\begin{equation*}
\aMt(n,2n,n) := \CM^{\tet^0}_\tau(n,2n,n) = \Rt(n,2n,n)/\!/_{\tet^0} \GL(n,2n,n).
\end{equation*}
The moduli space of Artin representations is highly non-reduced. In what follows, however, we will
only need a description of the underlying reduced scheme which we denote by $\aMt(n,2n,n)_\reduced$.
The proof of the following result is given in the appendix.

\prop{artin-moduli-k}
The map $W_\bu \mapsto \supp(\CH^3(W_\bu)))$ gives an $\SL(H)$-equivariant isomorphism
\begin{equation*}
\aMt(n,2n,n)_\reduced \cong S^n(\PP(H)).
\end{equation*}
\eprop

The dimension vector $(n,2n,n)$ that appears in the above proposition is a special case of the vector $\alpha(r,d,n)$ for $r = d = 0$.
More generally, for an arbitrary  triple $(r,d,n)$, we let
\begin{align*}
\UMt(r,d,n) := \CM^{\tet^0}_\tau(\alp(r,d,n)) = \Rt(\alp(r,d,n))/\!/_{\tet^0} \GL(\alp(r,d,n))
\end{align*}
be the moduli space of $\tet^0$-semistable $\alpha(r,d,n)$-dimensional representations
of $\bpt$, to be called the {\sf Uhlenbeck moduli space} of sheaves on $\BPt$
(the reason for this  name will become clear later).
%
%
%
We also consider the {\sf Gieseker moduli space} of sheaves on $\BPt$ defined as:
\begin{align*}
\GMt(r,d,n) := \Mt(\alp(r,d,n)) = \Rt(\alpha(r,d,n))/\!/_{\tnt}\GL(\alpha(r,d,n)),
\end{align*}
the moduli space of $\tnt$-semistable $\alpha(r,d,n)$-dimensions representations of  $\bpt$.
The reason for the name `Gieseker moduli' is motivated by the following result.

\prop{nevins}
The Gieseker moduli space $\GMt(r,d,n)$ is isomorphic to the moduli space of Gieseker
semistable sheaves on $\BPt$ constructed in~\cite{ns}. Moreover, the open subset
of $\GMt(r,d,n)$ of $\tnt$-stable representations corresponds, via the
isomorphism,
to the open set of
Gieseker stable sheaves    on $\BPt$.
\eprop
\prf
This follows immediately from~\refl{kronecker-stability}, as the functor of $\tnt$-(semi)stable representations of the quiver $\bpt$
is isomorphic to the functor of (semi)stable Kronecker complexes considered in~\cite{ns}.
\epr

\cor{sm}
If $\mathsf{gcd}(r,d,n)=1$, then $\GMt(r,d,n)$ is a fine moduli space;
moreover,
this moduli space is smooth.
\ecor
\prf
By~\cite[Prop.~7.15]{ns} the moduli space of semistable Kronecker complexes is fine. As the functor of semistable Kronecker complexes
is isomorphic to the functor of $\tnt$-semistable representations of $\bpt$, we conclude that $\GMt(r,d,n)$ is also a fine moduli space.
Moreover, the condition $\mathsf{gcd}(r,d,n) = 1$ implies that any $\tnt$-semistable representation of the quiver is $\tnt$-stable;
hence any Gieseker semistable sheaf is stable, and so the smoothness of the moduli space follows from~\cite[Thm.~8.1]{ns}.
\epr

%



\ssec{stratifications}{Stratifications}

Recall that by~\refp{tnt-ss} any $\tnt$-semistable representation $V_\bu$ can be written as $V_\bu(E)$
for a Gieseker semistable sheaf $E$.
This gives the decomposition of the moduli space $\GMt(r,d,n)$ corresponding to the length of $E^{**}/E$.
It will be shown in the appendix at the end of the paper that
this  decomposition is, in fact, an algebraic stratification. 

The proofs of   other results  of this subsection stated below
are also deferred to  Section 4.

\lem{stratification}
The Gieseker moduli space $\GMt(r,d,n)$ is naturally stratified by locally closed $\SL(H)$-invariant subsets
\begin{equation*}
\GMt(r,d,n) = \bigsqcup_{0\le k\le n} \GMt^k(r,d,n),
\end{equation*}
where the stratum $\GMt^k(r,d,n)$
corresponds to the locus of Gieseker semistable sheaves $E$ on $\BPt$ with $c_2(E^{**}/E)=k$.
\elem

The open stratum $\GMt^0(r,d,n) \subset \GMt(r,d,n)$ parameterizes locally free Gieseker semistable sheaves.

There is also an analogous stratification of the Uhlenbeck moduli space.

\lem{stratification-uhlenbeck}
The Uhlenbeck moduli space $\UMt(r,d,n)$ is naturally stratified by locally closed $\SL(H)$-invariant subsets
\begin{equation*}
\UMt(r,d,n) = \bigsqcup_{0\le k\le n} \UMt^k(r,d,n),
\end{equation*}
where the stratum $\UMt^k(r,d,n)$
corresponds to the locus of Mumford semistable sheaves $E$ on $\BPt$ with $c_2(E^{**}/E)=k$.
\elem

It should be pointed out that the strata of the {\em natural} stratifications of the Gieseker and Uhlenbeck moduli spaces \
are typically non-reduced. The reason for this  is that neither the scheme of  ``commutative points'' of $\BPt$
nor the moduli space of Artin sheaves are  reduced.
In the present paper we are only concerned with constructible sheaves on the schemes in question.
The difference between a scheme and the corresponding  reduced
scheme is therefore irrelevant. 
Thus, from now on we will abuse notation and write
$\GMt^k(r,d,n)$, resp. $\UMt^k(r,d,n)$, for the stratum of the
corresponding stratification equipped with the {\textsf{reduced}} scheme structure.

Observe that $\tet^0$-semistability is a degeneration  of $\tnt$-semistability.
Therefore, by standard results of geometric invariant theory
(see~\cite{DH}), there is a canonical $\SL(H)$-equivariant 
projective  morphism
\begin{equation*}
\gam_\tau\colon\GMt(r,d,n)\to\UMt(r,d,n).
\end{equation*}

\begin{Rem} We do not know how to define the morphism $\gam_\tau$ in
    terms
of coherent sheaves on $\BPt$, without using identifications of 
moduli spaces of coherent sheaves with the corresponding moduli spaces
of
quiver representations.
\end{Rem}

The main result of this section establishes  compatibility between the constructed statifications
of the Gieseker and Uhlenbeck moduli spaces and describes the relation between the strata.


\th{str}
(1) The map $\gamma_\tau:\GMt(r,d,n) \to \UMt(r,d,n)$ is compatible with the stratifications, i.e.,
$\gamma_\tau(\GMt^k(r,d,n)) \subset \UMt^k(r,d,n)$.

In the case where  $\tau\neq 0$ and the integers $r$ and $d$ are coprime, 
the  Gieseker compactification is smooth and the following holds:
\vskip 2pt 

(2) The  open set $\GMt^0(r,d,n)$ is the locus of  Gieseker stable
supermonadic representations; furthermore, this  open set
corresponds, via the isomorphism of Proposition \ref{P:nevins},
to the
locus of locally free Gieseker stable sheaves on $\BPt$.
Moreover, the
map  $\gamma_\tau$ yields an isomorphism $\GMt^0(r,d,n)\iso
\UMt^0(r,d,n)$.


(3) For any
$k > 0$ one has  an $\SL(H)$-equivariant isomorphism
\begin{equation*}
\UMt^k(r,d,n) \cong \GMt^0(r,d,n-k) \times \aMt(k,2k,k)_\reduced \cong \GMt^0(r,d,n-k) \times S^k\BP(H).
\end{equation*}
Using this isomorphism, for  $E\in
\GMt^k(r,d,n)$ we have
\begin{equation*}
\gamma_\tau(E) = (E^{**},\supp(E^{**}/E)).
\end{equation*}
In particular, 
the fiber of $\gamma_\tau$ over a point $(\CE,D) \in \GMt^0(r,d,n-k) \times S^k\BP^1 \subset \UMt(r,d,n)$
is the underlying reduced scheme for the moduli space of subsheaves $E \subset \CE$ with $\supp(\CE/E) = D$.
\eth

\begin{Rem} The relation between the moduli spaces $\GMt(r,d,n)$ and $\UMt(r,d,n)$
is analogous to that of the Gieseker and Uhlenbeck compactifications of the moduli
space of vector bundles on a commutative algebraic surface (this justifies our terminology).
There is, however, an important difference between  
commutative and noncommutative settings.
In the case of a commutative surface, each stratum of the  stratification of the  Uhlenbeck space
is isomorphic to a cartesian product of an  Uhlenbeck space of lower dimension and a symmetric powers of the surface.
In the case of a noncommutative surface, the corresponding stratum  of the  stratification of the  Uhlenbeck space
is isomorphic to a cartesian product of an  Uhlenbeck space of lower dimension and a 
 a symmetric power of a {\em curve}, rather than a surface.
This phenomenon will play a crucial role in subsequent sections.
\end{Rem}

\sec{ch6}{Rank $1$ sheaves and the Calogero-Moser space}

In this section we study the Gieseker and the Uhlenbeck moduli spaces of rank 1 and degree 0
torsion free sheaves on $\BPt$.

\ssec{calmos}{The compactifications}

To unburden the notation we write
$$
\obMt^n=\GMt(1,0,n),\quad \hbMt^n=\UMt(1,0,n), \quad
\obMt^{n,k}=\GMt^k(1,0,n), \quad
\hbMt^{n,k}=\UMt^k(1,0,n).
$$
It is well-known, cf. \cite[Prop.~8.13]{ns}, that the open strata $\GMt^{n,0} \subset \GMt^n$ and
$\UMt^{n,0} \subset \UMt$ can be identified with the Calogero-Moser
space.
Thus, the varieties $\obMt^n$ and $\hbMt^n$ provide
two {\em different}  compactifications  of the Calogero-Moser
space,
to be called the 
{\sf Gieseker and  Uhlenbeck compactifications}, respectively.
Furthermore, the variety $\obMt^n$ being smooth, the morphism $\gamma_\tau \colon \obMt^n \to \hbMt^n$ is a resolution of singularities.
Later on, we will use this resolution to compute the stalks of the IC sheaf of the Uhlenbeck compactification.


\th{isom}
For $\tau \ne 0$ we have $\SL(H)$-equivariant isomorphisms
\begin{equation*}
\GMt^{n,0} = \UMt^{n,0} \cong \sM_\tau^n.
\end{equation*}
\eth
\prf
The first isomorphism follows from~\refp{nevins} and~\cite[Prop.~8.13]{ns}. The second
is a consequence of the first and~\reft{str}(2).
\epr

In \refss{moduli-spaces}, we have introduced a contraction map $\gamma_\tau\colon\obMt^n\to\hbMt^n$.
By~\reft{str} it sends a torsion free sheaf $E$ to $(E^{**},\supp(E^{**}/E))$,
and for any $0\leq k\leq n$, we have
\begin{equation*}
\gamma_\tau(\obMt^{n,k})\ =\ \hbMt^{n,k}\ =\ \sM_\tau^{n-k}\x S^k\BP(H).
\end{equation*}
Below, we are going to describe the fibers of the map  $\gamma_\tau$.
%
%
Choose $0 \ne h \in H$ and let
\begin{equation*}
\BA^1_h := \BP(H) \smallsetminus \{h\}.
\end{equation*}
Thus, $\BA^1_h$ is an affine line and, for any $n\geq 0$, the set
$S^n\BA^1_h$ is  Zariski open and dense in $S^n\BP(H)$.
It is clear that these sets for all $h \in H$ form an open covering of $S^n\PP(H)$.

The zeroth Calogero-Moser space $\sM_\tau^0$ is, clearly,  a point.
This point corresponds,  under the isomorphism of~\reft{isom}, to the trivial line bundle $\CO_{\BPt}$.
For each nonzero vector $h \in H$, we consider the open subset
$\{\CO\} \times S^n\BA^1_h = \sM^0_\tau \times S^n\BA^1_h \subset \sM^0_\tau \times S^n\PP(H) = \hbMt^{n,n}$
and its preimage under the map $\gamma_\tau:\obMt^{n,n} \to \hbMt^{n,n}$:
\begin{equation}\label{equation-bnh}
B^n_h =
\gamma_\tau^{-1}(\{\calO\}\x S^n\BA^1_h).
\end{equation}
Analogously, we can take arbitrary locally free sheaf $\CE$ of rank 1 and degree 0, consider the locally closed subset
$\{\CE\} \times S^k\BA^1_h \subset \sM^m_\tau \times S^k\BA^1_h \subset \sM^m_\tau \times S^k\PP(H) = \hbMt^{m+k,k}$
and its preimage under the map $\gamma_\tau:\obMt^{m+k,k} \to \hbMt^{m+k,k}$:

\prop{dep}
For any locally free sheaf $\calE$ of rank $1$, $\calE\in\sM_\tau^{m}$,
there is an isomorphism
\begin{equation*}
\gamma_\tau^{-1}(\{\calE\}\x S^k\BA^1_h)\cong B^k_h.
\end{equation*}
\eprop
\prf
There is an integer $p \in \ZZ$ and two maps
\begin{equation*}
\phi:\CO(-p) \to \CE
\qquad\text{and}\qquad
\phi':\CE \to \CO(p)
\end{equation*}
such that $i^*(\phi) = h^p$ and $i^*(\phi') = h^p$. Indeed, take $p$ sufficiently large to have
\begin{equation*}
\Ext^1(\CE,\CO(p-1)) = \Ext^1(\CO(-p),\CE(-1)) = 0
\end{equation*}
and define $\phi$ and $\phi'$ as lifts of the compositions in the next two diagrams:
\begin{equation*}
\xymatrix@C=1em{
&& \CO(-p) \ar[r] \ar@{..>}[d]_\phi & i_*i^*\CO(-p) \ar[d]^{h^p} \\
0 \ar[r] & \CE(-1) \ar[r] & \CE \ar[r] & i_*i^*\CE \ar[r] & 0
}
\qquad
\xymatrix@C=1.5em{
&& \CE \ar[r] \ar@{..>}[d]_\phi' & i_*i^*\CE \ar[d]^{h^p} \\
0 \ar[r] & \CO(p-1) \ar[r] & \CO(p) \ar[r] & i_*i^*\CO(p) \ar[r] & 0
}
\end{equation*}

The maps $\phi$ and $\phi'$ give morphisms between the moduli spaces of surjections
$\CE \twoheadrightarrow F$ and $\CO \twoheadrightarrow F$ onto Artin sheaves $F$
of length $k$ with $\supp(F) \subset \BA^1_h$:
\begin{equation*}
(f:\CE \twoheadrightarrow F) \mapsto (f \circ \phi(p):\CO \twoheadrightarrow F)
\qquad\text{and}\qquad
(g:\CO \twoheadrightarrow F) \mapsto (g(p) \circ \phi':\CE \twoheadrightarrow F).
\end{equation*}
%
Here, in both cases we have identified $F$ with $F(p)$ via $h^p$. It is straightforward to check
that these maps are mutually inverse. As the preimages $\gamma_\tau^{-1}(\CE\times S^k\BA^1_h)$ and
$B^k_h = \gamma_\tau^{-1}(\CO\times S^k\BA^1_h)$ are identified by~\reft{str}(3) with the reduced schemes
underlying these moduli spaces, the constructed maps provide an isomorphism between them as well.
\epr

\newcommand{\disj}{{{\on{disj}}}}

The spaces $B^n_h$ come with a natural map $\gamma_\tau:B^n_h \to S^n\BA^1_h$.
In fact they enjoy the following factorization property. Define
the open subset $(S^{k_1}\BA^1\times S^{k_2}\BA^1)_\disj \subset S^{k_1}\BA^1\times S^{k_2}\BA^1$ as
\begin{equation}\label{disjoint-divisors}
(S^{k_1}\BA^1\times S^{k_2}\BA^1)_\disj = \{ (D_1,D_2) \in S^{k_1}\BA^1 \times S^{k_2}\BA^1 \ |\ \supp(D_1) \cap \supp(D_2) = \emptyset \}
\end{equation}
and
\begin{equation*}
(B_h^{k_1} \times B_h^{k_2})_\disj := (\gamma_\tau \times \gamma_\tau)^{-1}((S^{k_1}\AA_h^1 \times S^{k_2}\AA_h^1)_\disj) \subset B_h^{k_1} \times B_h^{k_2}.
\end{equation*}

\prop{factorization-b}
The collection of spaces $B^n_h$ has a factorization property, i.e., there is a collection of maps
\begin{equation*}
\psi_{k_1,k_2}:(B_h^{k_1}\times B_h^{k_2})_\disj \to B_h^{k_1+k_2}
\end{equation*}
for all positive integers $k_1,k_2$ which has the following properties: 
\begin{itemize}
\item (associativity) $\psi_{k_1+k_2,k_3} \circ (\psi_{k_1,k_2} \times \id) = \psi_{k_1,k_2+k_3} \circ (\id \times \psi_{k_1,k_2})$ for all $k_1,k_2,k_3$;
\item (commutativity) the maps $\psi_{k,k} \colon (B_h^k\times B_h^k)_\disj \to B_h^{2k}$ commute with the transposition of the factors on the source for all $k$;
\item (compatibility with the addition) the following diagram is Cartesian
\begin{equation*}
\xymatrix{
(B_h^{k_1}\times B_h^{k_2})_\disj \ar[rr]^-{\psi_{k_1,k_2}} \ar[d]_{\gamma_\tau\times\gamma_\tau} && B_h^{k_1+k_2} \ar[d]^{\gamma_\tau} \\
(S^{k_1}\BA^1\times S^{k_2}\BA^1)_\disj \ar[rr]^-{a_{k_1,k_2}} 	&& S^{k_1+k_2}\BA^1
}
\end{equation*}
where the bottom arrow is the addition morphism: $(D_1,D_2) \mapsto D_1 + D_2$.
\end{itemize}

\eprop
\prf
A point of $(B^{k_1}_h \times B^{k_2}_h)_\disj$ can be represented by a pair of Artin sheaves $F_1,F_2$ of length $k_1$ and $k_2$, respectively,
equipped with epimorphisms $\CO \twoheadrightarrow F_1$ and $\CO \twoheadrightarrow F_2$. 
Consider the sum $F = F_1 \oplus F_2$ and the map $\CO \to F$ given by the sum of the two above maps. Let us show it is surjective.
By~\refp{morph2}(2) it is enough to check that the map $\CO_{\PP(H)} \to i^*F = i^*F_1 \oplus i^*F_2$ is surjective.
The supports of the sheaves $i^*F_1$ and $i^*F_2$ being disjoint, this is equivalent to the surjectivity
of each of the maps $\CO_{\PP(H)} \to i^*F_1$ and $\CO_{\PP(H)} \to i^*F_2$. The latter holds by~\refp{morph2}(2).
Therefore, the sheaf $F$ equipped with the constructed epimorphism $\CO \twoheadrightarrow F$ gives a point of $B^{k_1+k_2}_h$.
Thus there is a well-defined morphism
\begin{equation*}
\psi^B_{k_1,k_2} \colon (B^{k_1}_h \times B^{k_2}_h)_\disj \to B^{k_1+k_2}_h.
\end{equation*}
We claim that it is a factorization isomorphism. Indeed, the associativity and the commutativity properties
are clear. It remains to check the compatibility with the addition operation, i.e., to check that the corresponding diagram is Cartesian.
The commutativity of the diagram follows from~\refl{artin-support}, so it remains to note that if $F$ is an Artin sheaf 
of length $k_1+k_2$ such that $\supp(F) = D_1 + D_2$ with disjoint divisors $D_1 \in S^{k_1}\AA^1_h$ and $D_2 \in S^{k_2}\AA^1_h$, 
then $F$ has a unique representation as a direct sum $F = F_1 \oplus F_2$ with $\supp(F_1) = D_1$
and $\supp(F_2) = D_2$ (this follows easily from~\refp{Art0}). 
%
%
%
%
\epr

The variety $B_h^k$ has a nice linear algebra description.
Fix a vector space $V$ of dimension~$k$. Let
\eq{rel}
\widetilde{B}{}^k_h = \{ (Y,Z,v)\in \End(V)\times \End(V)\times V\  |\  [Y,Z]=\tau Z^3\ \text{and $v$ is cyclic} \}.
\eeq
Here we say that a vector $v$ is {\sf cyclic} for a pair of matrices $(Y,Z)$ if
there is no proper vector subspace $V' \subset V$ that contains $v$ and is both $Y$-stable and $Z$-stable.
One has a natural  $\GL(V)$-action on $\widetilde{B}{}^k_h$
given by $g:\ (Y,Z,v)\mapsto (gYg^{-1}, gZg^{-1}, gv)$.

\th{26}
Th eaction of $\GL(V)$ on $\widetilde{B}{}^k_h$ is free and
\begin{equation*}
B^k_h \cong (\widetilde{B}{}^k_h/\GL(V))_\reduced.
\end{equation*}
Under this isomorphism the map $\gamma_\tau:B^k_h \to S^k\BA^1_h$ is induced by the map $\widetilde{B}^k_h \to S^k\BA^1_h$
which takes $(Y,Z,v)$ to $\Spec(Y)$.
\eth

\prf
Assume that $g\in \GL(V)$ acts trivially on a triple $(Y,Z,v)$. Let
$V^g\subset V$ be the space of invariants of $g$. Then $v\in V^g$ and
$Y(V^g)\subset V^g,\ Z(V^g)\subset V^g$, hence $V^g=V$ as $v$ is cyclic, and so $g=1$.


Consider the moduli space of surjections $\CO \twoheadrightarrow F$ such that $F$ an Artin sheaf
of length $k$ and $\supp(F) \subset \BA^1_h$. Let us show this moduli space is isomorphic to 
$\widetilde{B}{}^k_h/\GL(V)$. The theorem would then follow by  passing to the underlying reduced schemes.


Choose symplectic coordinates $x,y$ in $H$ such that the point $h \in \PP(H)$ is given by the equation $x = 0$.
Let $(Y,Z,v)$ be a point of $\widetilde{B}{}^k_h$.
Consider a graded vector space $V[x] := V \otimes \BC[x]$ with $\deg x = 1$, with $x$ acting by multiplications, and with the action of $y$ and $z$ defined by
\begin{equation*}
y = xY -\tau x^2Z^2\partial_x,
\qquad
z = xZ.
\end{equation*}
The commutation $[x,z] = 0$ is clear. Moreover, we have
\begin{equation*}
[y,z] = [xY -\tau x^2Z^2\partial_x, xZ] = x^2[Y,Z] - \tau x^2Z^2[\partial_x,x]Z = \tau x^2Z^3 - \tau x^2Z^3 = 0
\end{equation*}
and
\begin{equation*}
[x,y] = [x, xY -\tau x^2Z^2\partial_x] = - \tau x^2Z^2[x,\partial_x] = \tau x^2Z^2 = \tau z^2.
\end{equation*}
This shows that $V[x]$ is a graded $A^\tau$-module. Let $F$ be the corresponding coherent sheaf on $\BPt$.
By definition the map $x:V[x] \to V[x]$ is injective with finite-dimensional cokernel, hence the map $x:F \to F(1)$
is an isomorphism. In particular, the Hilbert polynomial $h_F(t)$ is constant, hence $F$ is an Artin sheaf
with $\supp(F) \subset \BA^1_h$ by~\refl{artin-support}. Moreover, the length of $F$ is equal to $\dim V = k$
and the vector $v \in V \subset V[x]$ gives a morphism  $A^\tau \to V[x]$ of $A^\tau$-modules. 
The cyclicity
assumption insures that this the map is surjective for all components of sufficiently large degree.
It follows that  the corresponding
morphism  $\CO \to F$ is surjective.
Note that the construction is $\GL(V)$-invariant.

Conversely, let $\CO \twoheadrightarrow F$ be a surjection where $F$ an Artin sheaf of length $k$
and $\supp(F) \subset \BA^1_h$. Choose an isomorphism $V \cong H^0(\BPt,F)$. 
By~\refl{artin-support}, the map $x:F \to F(1)$ is an isomorphism. Hence it 
induces an isomorphism $x \colon H^0(\BPt,F) \xrightarrow{\ \sim\ } H^0(\BPt,F(1))$ of the spaces of global sections.
Let $x^{-1}$ be the inverse isomorphism. We view $Y := x^{-1}y$ and $Z := x^{-1}z$, 
as endomorphisms of $V = H^0(\BPt,F)$. Finally, let $v$ be the image
of $1 \in H^0(\BPt,\CO)$ in $V$ under the map $\CO \twoheadrightarrow F$.

We claim that~\refe{rel} holds. First, we have $y = xY$, $z = xZ$ which gives relations
\begin{equation*}
x^2Z = xZx,
\qquad
xYxZ = xZxY,
\qquad\text{and}\qquad
x^2Y - xYx = \tau xZxZ.
\end{equation*}
It follows that
\begin{equation*}
x^2(YZ - ZY) =
x^2YZ - x^2ZY =
xYxZ + \tau xZxZ^2 - xZxY =
\tau x^2Z^3.
\end{equation*}
In the second equality we have used the third relation,
and in the third equality we have used the first two relations.
Since $x$ is an isomorphism, we deduce $[Y,Z] = \tau Z^3$. So, it remains to show that $v$ is cyclic. 
To this end, let  $V' \subset V$ be
an arbitrary subspace that contains $v$ and is closed under the action of $Y$ and $Z$.
The triple $(Y,Z,v)$, viewed as a triple for  the vector space $V'$, gives an Artin sheaf $F'$ and a surjection $\CO \twoheadrightarrow F'$.
The embedding $V' \subset V$ induces an embedding of sheaves $F' \hookrightarrow F$ and it is clear that the original map
$\CO \twoheadrightarrow F$ factors as $\CO \twoheadrightarrow F' \hookrightarrow F$. It follows that $F' = F$
and hence $V' = H^0(\BPt,F') = H^0(\BPt,F) = V$.

The two constructions are clearly mutually inverse and thus prove an isomorphism of moduli spaces
and hence the first part of the theorem.
For the second part it remains to show that $\supp(F) = \Spec(Y)$. But this is clear since
the action of the coordinate $y$ from $H^0(\BPt,F) = V$ to $H^0(\BPt,F(1)) = V$ is given by the operator $Y$.
\epr

We use the identification $B^k_h \cong\ (\widetilde{B}{}^k_h/\GL(V))_\reduced$
to investigate the properties of $B^k_h$.

\lem{nil}
If a pair $(Y,Z)$ satisfies \refe{rel} then $Z$ is nilpotent.
\elem
\prf
Note that $[Y,Z^p]=p\tau Z^{p+2}$ for $p \ge 3$. Since $\tau\ne0$, it follows that
$\Tr Z^{p+2}=0$ for any $p \ge 3$. Hence $Z$ is nilpotent.
\epr

\lem{ex}
For any nilpotent $Z$ there exist $Y$ and $v$ such that~\refe{rel} holds.
\elem
\prf
First, for any $u \in \BC$ take
\begin{equation*}
V = \C[t]/t^k,\qquad
Y = u + \tau t^3\partial_t,\quad
Z = t,\quad
v = 1.
\end{equation*}
Clearly, \refe{rel} holds,
so we have an example in case when $Z$ is just one Jordan block.
Note that $\Spec(Y) = ku \in S^k\BA^1_h$.

For arbitrary nilpotent $Z$ the Jordan decomposition of $Z$ is a direct sum
decomposition $Z = Z_1 \oplus \dots \oplus Z_m$ with blocks of size $k_1, \dots, k_m$.
Choosing $m$ distinct complex numbers $u_1,\dots,u_m$ we construct triples $(Y_i,Z_i,v_i)$
on $V_i = \BC[t]/t^{k_i}$ such that $\supp(Y_i,Z_i,v_i) = k_iu_i$. Factorization property of~\refp{factorization-b}
then shows that the direct sum $(\oplus Y_i,\oplus Z_i,\oplus v_i)$ is a point of $B^{k_1+\dots+k_m}_h$.
\epr


We will use a natural  one-to-one correspondence $\lam\mapsto O_\lam$,
between partitions of $k$ and the nilpotent conjugacy classes in
$\operatorname{End}(V)$, provided by Jordan normal form.
%
Let $B_h^\lam\subset B^k_h$ denote the set of all triples $(Y,Z,v)$ satisfying \refe{rel} with $Z\in O_\lam$.

\th{b-components}
We have a decomposition into a union of connected components
\begin{equation*}
B^k_h=\coprod_{\lam \in \frP(k)} B^\lam_h.
\end{equation*}
The component $B^\lam_h$ is smooth, connected and $k$-dimensional.
\eth
\prf
Let $O=O_\lam$ be a nilpotent orbit of the group $\GL(V)$. Let
\begin{equation*}
N^*_O\End(V)=\{(Y,Z)\ |\ [Y,Z]=0\}\subset\End(V) \times \End(V)
\end{equation*}
be the conormal bundle of $O$. Let
\begin{equation*}
_\tau\!N^*_O\End(V)=\{(Y,Z)\ |\ Z\in O,\ [Y,Z]=\tau Z^3\}.
\end{equation*}
Then according to~\refl{ex} the space $(_\tau\!N^*_O\End(V))_\reduced$ is a
$\GL(V)$-equivariant $N^*_O\End(V)$-torsor over~$O$.
In particular $(_\tau\!N^*_O\End(V))_\reduced$ is smooth and $k^2$-dimensional.

It is clear that for any pair $(Y,Z)$ the set of cyclic $v\in V$ is open.
Therefore,
\begin{equation*}
\widetilde{B}^k_h \subset \coprod_{\lam \in \frP(k)} \left(_\tau\!N^*_{O_\lam}\End(V)
\times V \right)
\end{equation*}
is an open subset.
Moreover, by~\refl{ex} it has a nonempty intersection with every component above.
The theorem now follows from~\reft{26}.
\epr

Consider the map $\gamma_\tau:B^k_h \to S^k\BA^1_h$.

\cor{dim}
The set of points $D \in S^k\BA^1_h$ such that $\dim \gamma_\tau^{-1}(D) \ge m$
has codimension at least $m$ in $S^k\BA^1_h$, and is empty for $m \ge k$.
\ecor
\prf
Since $B^k_h$ is equidimensional of dimension $k$ by~\reft{b-components}, it suffices to show that no component of $B^k_h$
is contained in the fiber of $\gamma_\tau:B^k_h \to S^k\BA^1_h$. For this note that the map $\gamma_\tau$
is equivariant with respect to the action of the group $\GG_a \subset \SL(H)$, the unipotent radical of the parabolic
which fixes $h \in H$, and that its action on $\BA^1_h$
is free.
%
\epr

\th{sml}
The map $\gamma_\tau\colon\GMt^n\to\UMt^n$ is small.
\eth
\prf
Let $(\UMt^n)_m \subset \UMt^n$ be the set of points over which the fiber of $\gamma_\tau$ has dimension~$m$.
Take any $0 \ne h \in H$. By~\refp{dep} for any $(\CE,D) \in \sM_\tau^{n-k} \times S^k\BA^1_h$ the fiber $\gamma_\tau^{-1}(\CE,D)$
is isomorphic to the fiber of the map $\gamma_\tau:B^k_h \to S^k\BA^1_h$ over $D$. In particular, by~\refc{dim}
the codimension of the set $(\UMt^n)_m \cap (\sM_\tau^{n-k} \times S^k\BA^1_h)$ in $\sM_\tau^{n-k} \times S^k\BA^1_h$
is at least $m$, and moreover $k > m$. Therefore
\begin{multline*}
\dim ((\UMt^n)_m \cap (\sM_\tau^{n-k} \times S^k\BA^1_h)) \le
\dim (\sM_\tau^{n-k} \times S^k\BA^1_h) - m \\
=2(n-k) + k - m = 2n - k - m < 2n - 2m.
\end{multline*}
Since the sets $\sM_\tau^{n-k} \times S^k\BA^1_h$ form an open covering of the stratum $\sM_\tau^{n-k} \times S^k\BA^1 = \UMt^{n-k,k}$
of a stratification of $\UMt^{n}$, the result follows.
%
\epr

\ssec{defo}{Deformation of $\GMt^n$ and $\UMt^n$}

The goal of this section is to show that the Gieseker and the Uhlenbeck compactifications
form a family over $\AA^1$ (with coordinate $\tau$) and check that the former is smooth.
%
To be more precise, consider the following graded algebra:
\begin{align*}
\sA =& \BC\langle x,y,z,\btau\rangle\Bigl.\Bigr/
\Bigl\langle[x,z]=[y,z]=[\btau,x]=[\btau,y]=[\btau,z]=0,[x,y]=\btau z ^2\Bigr\rangle, \\
& \deg x=\deg y=\deg z=1,\ \deg\btau=0.
\end{align*}
As $\btau$ is central of degree 0, this is an algebra over $\BC[\btau]$. In particular, we can specialize $\btau$
to any complex number $\tau$, which gives back the algebra $A^\tau$ we considered before.

Analogously, we consider the Koszul dual of $\sA$ over $\BC[\btau]$:
\begin{align*}
\sA^! =& \BC\langle \xi,\eta,\zeta,\btau \rangle/\langle \xi^2=\eta^2=\eta\xi+\xi\eta=\zeta\xi + \xi\zeta = \eta\zeta + \zeta\eta = \zeta^2 + \btau(\xi\eta - \eta\xi) = 0 \rangle,\\
& \deg \xi=\deg \eta = \deg \zeta = 1,\ \deg\btau=0.
\end{align*}
This is a graded $\BC[\btau]$-algebra. Note that each of its graded components $\sA^!_0$, $\sA^!_1$, $\sA^!_2$, $\sA^!_3$
is a free $\BC[\btau]$-module of finite rank (equal to $1$, $3$, $3$, and $1$ respectively).

%
Further, we consider the quiver $\bQ$ over $\BC[\btau]$ defined as
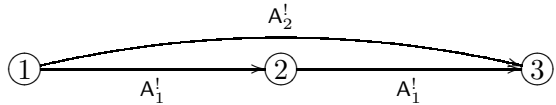
\begin{equation*}
\xymatrix{
*+[o][F]{1} \ar[rrr]_{\sA^!_1}\ar@/^1pc/[rrrrrr]^{\sA^!_2} &&&
*+[o][F]{2} \ar[rrr]_{\sA^!_1} &&&
*+[o][F]{3} }
\end{equation*}
(analogously to the quiver $\bpt$), and its representations in the category of $\BC[\btau]$-modules.
By definition such a representation is the data of three $\BC[\btau]$-modules $(\bV_1,\bV_2,\bV_3)$
and two morphisms of $\BC[\btau]$-modules $\bV_1 \otimes_{\BC[\btau]} \sA^!_1 \to \bV_2$ and
$\bV_2 \otimes_{\BC[\btau]} \sA^!_2 \to \bV_3$ such that the composition
$\bV_1 \otimes_{\BC[\btau]} \sA^!_1 \otimes_{\BC[\btau]} \sA^!_1 \to \bV_2 \otimes_{\BC[\btau]} \sA^!_1 \to \bV_3$
factors through $\bV_1 \otimes_{\BC[\btau]} \sA^!_2 \to \bV_3$.

Assuming each of $\bV_i$ is a free $\BC[\btau]$-module of finite rank, the space
\begin{equation*}
\Hom_{\BC[\btau]}(\bV_1 \otimes_{\BC[\btau]} \sA^!_1, \bV_2) \oplus \Hom_{\BC[\btau]}(\bV_2 \otimes_{\BC[\btau]} \sA^!_1, \bV_3)
\end{equation*}
is also a free $\BC[\btau]$-module of finite rank. We consider the associated vector bundle over $\Spec(\BC[\btau]) = \AA^1$
and its total space $\Tot(\Hom_{\BC[\btau]}(\bV_1 \otimes_{\BC[\btau]} \sA^!_1, \bV_2) \oplus \Hom_{\BC[\btau]}(\bV_2 \otimes_{\BC[\btau]} \sA^!_1, \bV_3))$
which is fibered over $\AA^1$ with fiber an affine space. The above factorization condition defines a Zarisky closed subspace
\begin{equation*}
\Rep_\bQ(\bV_\bu) \subset \Tot(\Hom_{\BC[\btau]}(\bV_1 \otimes_{\BC[\btau]} \sA^!_1, \bV_2) \oplus \Hom_{\BC[\btau]}(\bV_2 \otimes_{\BC[\btau]} \sA^!_1, \bV_3))
\end{equation*}
parameterizing all representations of the quiver $\bQ$ in $\bV_\bu$. By definition this is an affine variety over $\AA^1$.


%

Now we take a relatively prime triple $(r,d,n)$ such that~\eqref{assumptions-r-d-n} hold,
consider a triple of free $\BC[t]$-modules $(\bV_1,\bV_2,\bV_3)$ of ranks given
by the dimension vector $\alpha(r,d,n)$ of~\eqref{alpha}, and put
\begin{equation*}
\Rep_\bQ(\alp(r,d,n)) = \Rep_\bQ(\bV_\bu).
\end{equation*}
The group $\GL(\alp(r,d,n))$ acts naturally on the space $\Rep_\bQ(\alp(r,d,n))$
along the fibers of the projection $\Rep_\bQ(\alp(r,d,n)) \to \AA^1$.
Any rational polarization $\tet$ (in the sense of ~\refss{stabilo}) linearizes this action and thus
gives rise to the GIT quotient
\begin{equation*}
\CM^\tet_\bQ(\alpha(r,d,n)) = \Rep_\bQ(\alp(r,d,n)) /\!/_\tet \GL(\alp(r,d,n)).
\end{equation*}
By construction it comes with a map $\CM^\tet_\bQ(\alpha(r,d,n)) \to \AA^1$,
and clearly its fiber over a point $\tau \in \AA^1$ identifies with the moduli space
$\CM^\tet_\tau(\alpha(r,d,n))$.

Applying this construction in the case $\tet = \tet^0$, resp. $\tet = \tet^0 + \varepsilon\tet^1$, we construct
the following relative version of the Gieseker, resp. Uhlenbeck, compactification
\begin{equation*}
\GM(r,d,n) = \CM^{\tnt}_\bQ(\alpha(r,d,n))
\qquad\text{resp.}\qquad
\UM(r,d,n) = \CM^{\tet^0}_\bQ(\alpha(r,d,n)).
\end{equation*}
Standard results of GIT imply that there is  a morphism $\gamma:\GM(r,d,n) \to \UM(r,d,n)$
that commutes with the morphisms to $\AA^1$.

\prop{smd}
If $\operatorname{gcd}(r,d,n)=1$, then the map $\GM(r,d,n) \to \AA^1$ is smooth and projective.
In particular, $\GM(r,d,n)$ is a smooth variety.
\eprop
\prf
By~\refl{kronecker-stability} the moduli space $\GM(r,d,n)$ coincides with the moduli space of Gieseker semistable
sheaves of rank $r$, degree $d$ and second Chern class $n$ for the family $\sA$ of Artin--Schelter algebras over $\BC[\btau]$
constructed in~\cite{ns}. The smoothness and the projectivity of the latter is proved in Theorem~8.1 of {\em loc.\ cit.}
\epr

In  the special case $r = 1$, $d = 0$, we will use simplified notation $\GM^n=\GM(1,0,n)$, resp.  $\UM=\UM(1,0,n)$.

%
%

\ssec{fixed}{Fixed points}

We choose a torus $T \subset \SL(H)$ and consider its action on the Calogero-Moser space
and its Gieseker and Uhlenbeck compactifications. The stratifications and the map $\gam$
are $\SL(H)$-equivariant and hence $T$-equivariant as well. We aim at a description
of the set of $T$-fixed points on $\UMt^n$. Recall first what is known about the $T$-fixed locus of $\sM^n_\tau$.

\lem{fixed-cm}(\cite[Proposition~6.11]{W})
For $\tau \ne 0$ the set of $T$-fixed points in $\sM^n_\tau$ is in a natural bijection
with the set $\frP(n)$ of partitions of $n$.
\elem

We denote by $\bc^n_\lambda \in \sM^n_\tau$ the $T$-fixed point corresponding to a partition $\lambda \in \frP(n)$.
In particular, $\bc^0 \in \sM^0_\tau$ is the unique point (it is automatically $T$-fixed).
Let also $P_0,P_\infty \in \PP(H)$ be the $T$-fixed points on the line $\PP(H)$.

\lem{fixed-uhlenbeck}
For any $\tau \ne 0$ the set of $T$-fixed points in $\UMt^n$ is finite. Moreover,
\begin{equation*}
(\UMt^n)^T = \{ (\bc^m_\lambda, k_0 P_0 + k_\infty P_\infty) \in \sM^m_\tau \times S^{n-m}\BP^1 \mid \lambda\in\frP(m),\ k_0+k_\infty=n-m\}.
\end{equation*}
\elem
\prf
It is enough to describe $T$-fixed points on each of the strata $\sM^m_\tau \times S^{n-m}\PP(H)$ of the stratification of $\UMt^n$.
As the product decomposition is $\SL(H)$-invariant, it is enough to describe fixed points on each factor.
On first factor we use~\refl{fixed-cm}, and on $S^k\PP(H)$ a description of fixed points is evident.
\epr

Recall that a $T$-fixed point $P$ is called {\sf attracting} if all the weights of the $T$-action
on the tangent space at point $P$ are positive.

\lem{attracting-uhlenbeck}
The Uhlenbeck compactification $\UMt^n$ of the Calogero-Moser space has a unique attracting $T$-fixed
point $(\bc^0,nP_0) \in \sM^0_\tau\times S^n\PP(H) \subset \UMt^n$.
\elem
\prf
Since $\UMt^n$ is a projective variety, the $T$-action on it should have at least one attracting point.
On the other hand, $\sM^n_\tau$ is a sympletic manifold, and the $T$-action preserves
the sympletic structure~\cite{kks}, hence for $m > 0$ the weights of $T$ on the tangent
spaces at points $\bc^m_\lambda$ are pairwise opposite, and thus for $m > 0$ the $T$-fixed
points $(\bc^m_\lambda,k_0 P_0 + k_\infty P_\infty)$ are not attracting. Therefore, each
attracting point of the $T$-action on $\UMt^n$ lies on $\sM^0\times S^n\PP(H) = S^n\PP(H)$.
As it also should be an attracting point for the $T$-action on $S^n\PP(H)$, it should
coincide with $(\bc^0,nP_0)$.
\epr

\ssec{comp}{The IC sheaf of the Uhlenbeck compactification}


In this section we will prove~\reft{main}. 
The statement of the theorem and the arguments we use are purely topological.
We refer to~\cite{bbd} for the notion of IC sheaf and the general machinery.

We start with computing the stalks of the IC-sheaf at the deepest stratum of 
the Uhlenbeck stratification.
Since for $n=0$ the Calogero-Moser space $\sM^0_\tau$ is just a point, by~\reft{str}(2) we have
$S^n\PP(H) = S^n\PP(H) = \sM^0_\tau \times S^n\PP(H) \subset \UMt^n$.
Recall also the diagonal stratification~\eqref{diagonal-stratification} of $S^n\PP(H)$
and its deepest stratum $S_{(n)}\PP(H) \subset S^n\PP(H)$.




\prop{one stalk}
For any $P \in \PP(H)$ the stalk of the sheaf $\ic(\UMt^n)$ at the point $(\CO,nP)$ of the stratum $\sM^0_\tau \times S_{(n)}\PP(H) \subset \UMt^n$ is
isomorphic to
\eq{chow}
\ic(\UMt^n)_{(\CO,nP)} = \bigoplus_{\mu\in\frP(n)}\BC[2l(\mu)].
\end{equation}
\eprop
\prf
Let $T \subset \SL(H)$ be a torus such that $P = P_0$ is the attracting point for the action of~$T$ on $\PP(H)$.
The computation
is based on the following ``deformation diagram'':
\begin{equation*}
\xymatrix{
\GM^n_0\ar[d]^<>(0.5){\gamma_0}\ar@{^{(}->}[rr]^<>(0.5){\tilde\varsigma} &&
\GM^n\ar[d]^<>(0.5){\gamma} &&
\GM^n_\bseta\ar[d]^<>(0.5){\gamma_\bseta}\ar@{_{(}->}[ll]
\\
\UM^n_0\ar[d]\ar@{^{(}->}[rr]^<>(0.5){\varsigma} &&
\UM^n\ar[d]^<>(0.5){p} &&
\UM^n_\bseta\ar[d]\ar@{_{(}->}[ll]
\\
\{0\}\ar@{^{(}->}[rr] &&
\BA^1\ar@/^/[u]^<>(0.5){\sigma} &&
\BA^1\sminus\{0\}\ar@{_{(}->}[ll]
}
\end{equation*}
Here the middle column is the deformation family over $\AA^1$ of~\refp{smd}
with $p$ being the structure map. The left column is the fiber over the point $0 \in \AA^1$,
while the right column is the base change to $\AA^1 \setminus \{0\} \subset \AA^1$. Finally,
the map $\sigma:\AA^1 \to \UM^n$ is defined as follows.

For any $\tau \ne 0$ we put $\sigma(\tau) = n\cdot P \in \UMt^n \subset \UM^n$.
Clearly, this is a regular map $\AA^1 \setminus \{0\} \to \UMt$. 
Since $\UM^n$ is proper over $\AA^1$, the map extends to a map $\sigma:\AA^1 \to \UM^n$. 
By construction, the restriction of $\sigma$  to $\AA^1 \setminus\{0\}$ is a section of the map $p$.
It follows by continuity that $\sigma$  is a section of  $p$.



Let $\C^\times_{\leq1}$  be a sub-semigroup of $\BC^\times = T$
formed by the complex numbers with absolute value $\leq1$, and
let $F = \sigma(\BA^1) \sset \UM^n$ be the image  of the section $\sigma$.
Further, let $\CU\subset \UM^n$ be a small open neighborhood (in the analytic topology)
of the point $\sigma(0)$. Without loss of generality,
we may choose the set $\CU$ to be  $\C^\times_{\leq1}$-stable.
Note that $F$ is the attracting connected component
of $(\UM^n)^T$, by \refl{attracting-uhlenbeck}.
Therefore, shrinking  $\CU$ further, if necessary,
one may assume  in addition that we have $\CU^T=F\cap \CU$.
The action of $\C^\times_{\leq1}$ preserves the fibers of $p:\ \CU\to\BA^1$,
and contracts $\CU$ to the section $F=\sigma(\BA^1)$. According
to~\cite[Lemma~6]{Br},
for any $\C^\times$-equivariant complex $\CF$ of constructible sheaves on~$\UMt^n$,
the natural morphism $\sigma^*\CF\to p_*(\CF|_\CU)$ is an isomorphism.
In other words, for any $\tau\in\BA^1$, there is a natural isomorphism
\eq{is0}
H^\bullet(\CU\cap \UMt^n, \ \CF)\cong \CF|_{\sigma(\tau)}.
\eeq


Next, let $\psi_p$, resp. $\psi_\pg$,  denote the nearby cycles
functor~\cite[8.6]{KS}
with respect to the function $p$, resp. $\pg$. Note  that the morphism $\pg$
being smooth, we have
$\psi_\pg(\unl\C{}_{\GM^n_\bseta})=\unl\C{}_{\GM^n_0}$.
Therefore, using the proper base change for
nearby cycles (see e.g.,~\cite[Exercise~VIII.15]{KS}) we obtain
\begin{equation*}
(\gamma_0)_*\unl\C{}_{\GM^n_0}=(\gamma_0)_*(\psi_\pg(\unl\C{}_{\GM^n_\bseta}))
=\psi_p((\gamma_\bseta)_*\unl\C{}_{\GM^n_\bseta}).
\end{equation*}
The map  $\gamma_\bseta: \GM^n_\bseta\to \UM^n_\bseta$
 is  a small and proper morphism.
Hence, we have an isomorphism $\IC(\UM^n_\bseta) \cong (\gamma_\bseta)_*\unl\C{}_{\GM^n_\bseta}[2n]$.
Combining the above isomorphisms and
taking stalks at the point $\sigma(0)$ yields
\eq{is1}
\big((\gamma_0)_*\unl\C{}_{\GM_0^n}\big)|_{\sigma(0)}\cong
\big(\psi_p((\gamma_\bseta)_*\unl\C{}_{\GM^n_\bseta})\big)|_{\sigma(0)}
\cong\big(\psi_p(\IC(\UM^n_\bseta))\big)|_{\sigma(0)}[-2n].
\eeq

Further, by definition
of the functor $\psi_p$, for a sufficiently small open set $\CU$
as above and for any $\tau\neq 0$
with a sufficiently small absolute value,
one has $$
\big(\psi_p(\IC(\UM^n_\bseta))\big)|_{\sigma(0)}
\cong H^\bullet\big(\CU\cap p^{-1}(\tau),\ \IC(\UM^n_\bseta)\big)
\cong H^\bullet\big(\CU\cap\UMt^n,\ \IC(\UMt^n)\big).
$$

Thus, comparing the LHS and the RHS in \refe{is1},
we obtain
\begin{align}H^\bullet(\gamma_0^{-1}(\sigma(0)))[2n]\cong
\big((\gamma_0)_*\unl\C{}_{\GM_0^n}\big)|_{\sigma(0)}[2n]&\cong
\big(\psi_p(\IC(\UM^n_\bseta))\big)|_{\sigma(0)}\label{is4}\\
&\cong
H^\bullet\big(\CU\cap\UMt^n,\ \IC(\UMt^n)\big)
\cong\IC(\UMt^n)|_{\sigma(\tau)},\nonumber
\end{align}
where the last isomorphism is a special case of  \refe{is0}
for $\CF=\IC(\UMt^n)$.

To complete the proof, we observe that the
fiber $\gamma_0^{-1}(\sigma(0)))$ is the ``central fiber''
of the Hilbert--Chow morphism $\Hilb^n(\BP^2)\to S^n\BP^2$.
In other words, the variety $\gamma_0^{-1}(\sigma(0)))$
is nothing but $\Hilb^n_0(\BA^2)$, the punctual Hilbert scheme of
$n$ infinitesimally close points in $\BA^2$.
The Betti numbers of the punctual Hilbert scheme
are well-known, cf. e.g.,~\cite{N}. Specifically,
all odd Betti numbers vanish and one
has the formula
\[\dim H^{2k-2}(\Hilb^n_0(\BA^2))=\#\{\mu\in \frP(n)\mid l(\mu)=k\}.\]
It follows from (\ref{is4}) that, for
$\tau$ sufficiently small, the dimensions of the cohomology
groups of the stalk $\IC(\UMt^n)|_{\sigma(\tau)}$
are given by the same formula. This is equivalent to the
statement of the proposition.
\epr


Now~\reft{main} follows from~\refp{one stalk} and the factorization
property of~\refp{factorization-b}. In effect, due to $\SL(H)$-equivariance,
it suffices to find the stalks of $\IC(\UMt^n)$ at $S_\lambda\BA^1_h\subset S_\lambda\PP(H)$.
Given a point $\CE \in \sM^m$ and $D\in S_\lambda\BA^1_h$, due to the smallness of $\gamma_\tau$,
\begin{equation*}
\IC(\UMt^n)_{(\CE,D)} = H^\bullet(\gamma_\tau^{-1}(\{\CE\}\times D),\BC).
\end{equation*}
Further, by~\refp{dep} we have
\begin{equation*}
H^\bullet(\gamma_\tau^{-1}(\{\CE\}\times D),\BC) \cong H^\bullet(\gamma_\tau^{-1}(\{\CO\}\times D),\BC).
\end{equation*}
Now if $D = \sum k_iP_i$ with pairwise distinct points $P_i \in \BA^1_h$,
then according to~\refp{factorization-b},
$\gamma_\tau^{-1}(\{\CO\}\times D)\cong \prod_i\gamma_\tau^{-1}(\{\CO\}\times k_iP_i)$,
and hence
\begin{equation*}
H^\bullet(\gamma_\tau^{-1}(\{\CO\}\times D),\BC)\cong
\bigotimes_i H^\bullet(\gamma_\tau^{-1}(\{\CO\}\times k_iP_i),\BC).
\end{equation*}
Due to the smallness of $\gamma_\tau$,
$H^\bullet(\gamma_\tau^{-1}(\{\CO\}\times k_iP_i),\BC)=\IC(\UMt^i)_{(\CO,k_iP_i)}$,
and the latter stalk is known from~\refp{one stalk}. This completes the proof of~\reft{main}. \qed

\sec{app}{Appendix}
\noindent
In this Appendix we collect the proofs of some results from     ~\refs{ch5}.
Throughout, we assume that $\tau\neq 0$.
\ssec{artin-proofs}{Moduli spaces of Artin sheaves}

Let $\PP^1_3$ be the third infinitesimal
neighborhood of the line at infinity $\PP(H)$ in $\BPt$,
i.e., the projective spectrum of a commutative graded algebra $\C[x,y,z]/z^3$.


\lem{artin-moduli}
The moduli space $\aMt(1,2,1)$ is a fine moduli space.
It is isomorphic to the third infinitesimal neighborhood of a line on a plane:
$\aMt(1,2,1) \cong \PP^1_3$.
\elem
\prf
The data of a $(1,2,1)$-dimensional representation of $\bpt$ amounts to giving a pair of maps
\begin{equation*}
\C \xrightarrow{f} \C^2\otimes A^\tau_1
\qquad\text{and}\qquad
\C^2 \xrightarrow{g} \C\otimes A^\tau_1
\end{equation*}
such that the composite map 
\begin{equation*}
\C \xrightarrow{f} \C^2\otimes A^\tau_1 \xrightarrow{g\otimes 1} \C \otimes A^\tau_1\otimes A^\tau_1 \to \C \otimes A^\tau_2
\end{equation*}
is equal to zero.
Put 
$\sigma := (g\otimes 1)(f(1)) $ and let $K$ be the kernel of the multiplication
map $A^\tau_1\otimes A^\tau_1 \to A^\tau_2$.
Then, the last condition holds iff $\sigma$ belongs to $K$.

We have a natural identification $A^\tau_1\otimes A^\tau_1 = \Hom((A^\tau_1)^*,A^\tau_1)$,
so  one may view $\sigma$ as a map $(A^\tau_1)^*\to A^\tau_1$.
Also, one may view $f$ as a map $f^T\colon \BC^2 \to \BC\otimes A^\tau_1$.
Then, the  representation  of $\bpt$ associated with $(f,g)$ is $\tet^0$-semistable iff
each of the  maps $f^T$ and $g$ is injective, which holds iff  $\sigma$, viewed as a map,
 has rank~2
(note that the $\C^2$-component of the representation is just the image of $\sigma$). Thus the moduli space is nothing
but the degeneration scheme of the morphism
\begin{equation*}
(A^\tau_1)^* \otimes \CO_{\PP(K)}(-1) \to A^\tau_1 \otimes \CO_{\PP(K)}
\end{equation*}
on $\PP(K)$. As $K \subset A^\tau_1 \otimes A^\tau_1$ can be written as
\begin{equation*}
K = \{ u(y\otimes z - z\otimes y) + v(x\otimes z - z\otimes x) + w(x\otimes y - y\otimes x - \tau z\otimes z) \mid u,v,w \in \BC \},
\end{equation*}
the above morphism is given by the matrix
\begin{equation}\label{matrix-u-v-w}
\left(
\begin{matrix}
0 & w & v \\
-w & 0 & u \\
-v & -u & -\tau w
\end{matrix}
\right)
\end{equation}
and the degeneration condition is given by its determinant which is equal to
\begin{equation*}
\det\left(
\begin{matrix}
0 & w & v \\
-w & 0 & u \\
-v & -u & -\tau w
\end{matrix}
\right) = -\tau w^3.
\end{equation*}
This means that the moduli space is the subscheme of $\PP(K)$ given by the equation $w^3 = 0$, i.e.,
the third infinitesimal neighborhood $\PP^1_3$ of the line $\PP^1 = \{w = 0\}$ in the plane $\PP^2$.

To show that the moduli space is fine we should construct a universal family. For this we restrict the map
$(A^\tau_1)^* \otimes \CO_{\PP(K)}(-1) \to A^\tau_1 \otimes \CO_{\PP(K)}$
to $M := \PP^1_3$. This is a morphism of constant rank $2$ (the rank does not drop to 1 since among the 2-by-2 minors of the matrix~\eqref{matrix-u-v-w} one easily finds $u^2$, $v^2$, and $w^2$),
hence its image is a rank $2$ vector bundle~$\CV_2$. It comes equipped
with a surjective map $(A^\tau_1)^* \otimes \CO_{M}(-1) \to \CV_2$ and an injective map
$\CV_2 \to A^\tau_1\otimes \CO_M$. Clearly these two maps provide $(\CO_M(-1),\CV_2,\CO_M)$ with a structure of a family of representations of the quiver~$\bpt$.
The above arguments show it is a universal family.
\epr


%
%
%
%
%

Now we give a description of the reduced structure of the space $\aMt(k,2k,k)$ for $k > 1$.

\prf[Proof of~\refp{artin-moduli-k}]
Consider the subset $\aRt^{\tet^0}(k,2k,k) \subset \aRt(k,2k,k)$ of all $\tet^0$-semistable $(k,2k,k)$-dimensional representations of $\bpt$
and let $\CW_\bullet$ be the universal representation of the quiver over $\BPt$. Let $\CF$ be the universal
sheaf on the product $\aRt^{\tet^0}(k,2k,k) \times \BPt$, i.e., the sheaf defined by exact sequence
\begin{equation*}
0 \to \CW_1 \boxtimes \CO(-1) \to \CW_0 \boxtimes \CO \to \CW_1 \boxtimes \CO(1) \to \CF \to 0.
\end{equation*}
Then the support map defined in~\refl{artin-support} gives a map $\supp:\aRt^{\tet^0}(k,2k,k) \to S^k\PP^1$.
The map is clearly $\GL(k,2k,k)$-equivariant, hence descends to a map from the moduli space
\begin{equation*}
\supp:\aMt(k,2k,k) \to S^k\PP^1.
\end{equation*}
On the other hand, we clearly have an embedding which takes a $k$-tuple of $(1,2,1)$-dimensional Artin representations
$W^1_\bullet$, $W^2_\bullet$, \dots, $W^k_\bullet$ to their direct sum
\begin{equation*}
(\aRt(1,2,1))^k \to \aRt(k,2k,k),
\qquad
(W^1_\bullet, W^2_\bullet, \dots, W^k_\bullet) \mapsto W^1_\bullet \oplus W^2_\bullet \oplus \dots \oplus W^k_\bullet.
\end{equation*}
This map is equivariant with respect to the action of the group $\GL(1,2,1)^k \rtimes \fS_k$ on the source,
such that the $i$-th factor $\GL(1,2,1)$ acts naturally on the $i$-th factor of $(\aRt(1,2,1))^k$ and $\fS_k$ permutes the factors,
and the action on the target is given by
a natural embedding $\GL(1,2,1)^k \rtimes \fS_k \subset \GL(k,2k,k)$. The $\Proj$ construction of the GIT quotient implies that the map
induces a morphism of the GIT quotients
\begin{equation*}
(\aRt(1,2,1))^k /\!/_{\tet^0} (\GL(1,2,1)^k \rtimes \fS_k) \to \aRt(k,2k,k) /\!/_{\tet^0} \GL(k,2k,k).
\end{equation*}
The quotient on the right is just the moduli space $\aMt(k,2k,k)$. The quotient on the left
can be identified with $(\aMt(1,2,1))^k / \fS_k$,
so it is isomorphic to $S^k(\PP^1_3)$ by~\refl{artin-moduli}.
Restricting to the reduced subscheme, we obtain a map
\begin{equation*}
\Sigma \colon S^k\PP^1 = S^k(\PP^1_3)_\reduced \to \aMt(k,2k,k).
\end{equation*}
We are going to show that the constructed maps $\supp$ and $\Sigma$ induce isomorphisms between $S^k\PP^1$ and
the reduced moduli space $\aMt(k,2k,k)_{\mathrm{red}}$.

For this we note that the maps give bijections between the sets of closed points of $S^k\PP^1$ and $\aMt(k,2k,k)$,
since by~\refp{Art0}(3) any Artin sheaf is S-equivalent to a direct sum of structure sheaves for a unique collection
of points (which are given back by the support map). Note also that both $S^k\PP^1$ and $\aMt(k,2k,k)$ are projective
varieties, hence the map $\Sigma$ is proper. Finally, $S^k\PP^1 \cong \PP^k$ is normal.

So, it is enough to show that any proper regular map from a reduced normal scheme to a reduced scheme inducing
a bijection on the sets of closed points is an isomorphism. Locally, we just have an integral (due to properness)
extension of rings with the bottom ring being integrally closed (by normality), hence it is an isomorphism.
\epr


\newcommand{\GRt}{{}^G\Rt}
\newcommand{\URt}{{}^U\Rt}

\ssec{statifications-proofs}{Stratifications}

Here we construct the required stratifications of the Gieseker and Uhlenbeck moduli spaces.

\prf[Proof of~\refl{stratification}]
Let $\GRt := \osRt^{\tnt}(\alpha(r,d,n)) \subset \osRt(\alpha(r,d,n))$ be the open subset of
$\tnt$-semistable $\alpha(r,d,n)$-dimensional representations of $\bpt$.
Let $\CV_\bu$ be the universal family of representations over $\GRt$. Consider the universal monad
\begin{equation*}
\CV_1 \boxtimes \CO(-1) \to \CV_2 \boxtimes \CO \to \CV_3 \boxtimes \CO(1)
\end{equation*}
on $\GRt \times \BPt$ and denote its cohomology sheaf by $E$. For each point $s \in \GRt$
we denote by $E_s$ the restriction of $E$ to $\{s\} \times \BPt$. Note that this is just
the cohomology sheaf of the monad $\CV_{1s}\otimes \CO(-1) \to \CV_{2s} \otimes \CO \to \CV_{3s}\otimes\CO(1)$.
In particular, the sheaf $E$ is flat over $\GRt$.
%
%

Consider also the dual monad on $S\times\BPt$
\begin{equation*}
\CV_3^\vee\boxtimes\CO(-1) \to \CV_2^\vee \boxtimes\CO \to \CV_1^\vee \boxtimes\CO(1)
\end{equation*}
and let $\CF$ be the cokernel of the last map
\begin{equation*}
\CF := \Coker( \CV_2^\vee \boxtimes\CO \to \CV_1^\vee \boxtimes\CO(1) ).
\end{equation*}
For each point $s \in S$ we have
\begin{equation*}
\CF_s \cong \Coker( \CV_{2s}^\vee \boxtimes\CO \to \CV_{1s}^\vee \boxtimes\CO(1) ) \cong \underline{\Ext}^1(E_s,\CO) \cong \underline{\Ext}^2(E_s^{**}/E_s,\CO).
\end{equation*}
Thus it is an Artin sheaf, but its length may vary from point to point.
Consider the flattening stratification of $S$ for $\CF$:
\begin{equation*}
\GRt = \GRt^{\ge 0} \supset \GRt^{\ge 1} \supset \GRt^{\ge 2} \supset \dots \supset \GRt^{\ge n} \supset \GRt^{\ge n+1} = \emptyset,
\end{equation*}
where $\GRt^{\ge k}$ is the subscheme of points $s \in \GRt$ where the length of $\CF_s$ is at least $k$.
This stratification is $\GL(\alpha(r,d,n))$-invariant, so it gives a stratification
of the GIT quotient $\GRt/\!/_{\tnt}\GL(\alpha(r,d,n))$, i.e., of the Gieseker moduli space $\GMt(r,d,n)$.
Finally, we replace each stratum by its underlying reduced subscheme.
%
%
\epr
%


Below we will need the following result on universal families.

\prop{family-on-strata}
Let $\CV_\bu$ be the universal family of $\tnt$-semistable $\alpha(r,d,n)$-dimensional representations
of $\bpt$ over $\GRt^k := \GRt^{\ge k} \setminus \GRt^{\ge k+1}$.
Then there is a natural exact sequence
\begin{equation*}
0 \to \CW_\bu \to \CV_\bu \to \CU_\bu \to 0
\end{equation*}
of families of representations over $\GRt^k$ where $\CW_\bu$ is a family of Artin representations
of dimension $(k,2k,k)$ and $\CU_\bu$ is a family of supermonadic $\tnt$-semistable representations.
\eprop
\prf
We freely use the notation introduced in the proof of~\refl{stratification}.
By assumption, $\CF$ is a flat (over~$\GRt^k$)
family of Artin sheaves of length $k$.
Let
\begin{equation*}
\CW'_1 \boxtimes \CO(-1) \to \CW'_2 \boxtimes \CO \to \CW'_3 \boxtimes \CO(1)
\end{equation*}
be the Beilinson resolution of $\CF$.
The family $\CF$ being flat, it follows that each of the sheaves
$\CW'_1$, $\CW'_2$, and $\CW'_3$ is a vector bundle of rank $k$, $2k$, and $k$, respectively.
The functoriality of the Beilinson resolution yields a morphism of resolutions
\begin{equation*}
\xymatrix{
\CV_3^\vee \boxtimes\CO(-1) \ar[r] \ar[d] & \CV_2^\vee \boxtimes\CO \ar[r] \ar[d] & \CV_1^\vee \boxtimes\CO(1) \ar[d] \\
\CW'_1 \boxtimes\CO(-1) \ar[r] & \CW'_2 \boxtimes\CO \ar[r] & \CW'_3 \boxtimes\CO(1).
}
\end{equation*}

We claim that the induced morphisms $\CV_i^\vee \to \CW'_{4-i}$, of vector bundles on $\GRt^k$, are surjective.
It suffices to prove this pointwise, i.e., for a single representation rather than a family. In that case both
$\CV^\vee_\bu$ and $\CW'_\bu$ are $\tet^0$-semistable. Hence so is the image of the map and
also the  corresponding sheaf $\CF' \subset \CF$. A $\tet^0$-semistable
subrepresentation of an Artin representation is Artin.
Thus, $\CF'$ is an Artin sheaf such that the map $\CV^\vee_1 \to \CF$ factors through  $\CF'$. 
This forces $\CF' = \CF$,  by definition of $\CF$, and our claim follows.


Let $\CW_i = (\CW'_{4-i})^\vee$ and $\CU_i = \Ker(\CV_i^\vee \to \CW'_{4-i})^\vee$, so that we have an exact sequence of monads
\begin{equation*}
\xymatrix{
\CW_1 \boxtimes\CO(-1) \ar[r] \ar[d] & \CW_2 \boxtimes\CO \ar[r] \ar[d] & \CW_3 \boxtimes\CO(1) \ar[d] \\
\CV_1 \boxtimes\CO(-1) \ar[r] \ar[d] & \CV_2 \boxtimes\CO \ar[r] \ar[d] & \CV_3 \boxtimes\CO(1) \ar[d] \\
\CU_1 \boxtimes\CO(-1) \ar[r] & \CU_2 \boxtimes\CO \ar[r] & \CU_3 \boxtimes\CO(1).
}
\end{equation*}
By construction, the family $\CU_\bu$ is $\tnt$-semistable and supermonadic, so this exact sequence
is the one we need.
\epr

We are ready to construct the stratification of the Uhlenbeck moduli space.
The situation here is a bit more complicated than in the Gieseker case, since Artin
representations can appear both as subrepresentations and as quotient representations
of a $\tet^0$-semistable representation. We 
deal with the latter complication first.

\prf[Proof of~\refl{stratification-uhlenbeck}]
Let $\URt := \hsRt^{\tet^0}(\alpha(r,d,n)) \subset \hsRt(\alpha(r,d,n))$ be the open subset of
$\tet^0$-semistable $\alpha(r,d,n)$-dimensional representations of $\bpt$. Let $\CV_\bu$ be
the universal family of representations over $\URt$.
Consider the family of sheaves $\CF' := \Coker(\CV_2 \boxtimes \CO \to \CV_3 \boxtimes \CO(1))$ over $\URt \times \BPt$.
Note that these are Artin sheaves of length at most $n$. Let
\begin{equation*}
\URt = \URt^{\ge 0,\bu} \supset \URt^{\ge 1,\bu} \supset \URt^{\ge 2,\bu} \supset \dots \supset \URt^{\ge n,\bu} \supset \URt^{\ge n+1,\bu} = \emptyset,
\end{equation*}
be the flattening stratification for the sheaf $\CF'$. 
We restrict the family $\CV_\bu$ to each stratum $\URt^{k,\bu} = \URt^{\ge k,\bu} \setminus \URt^{\ge k+1,\bu}$
and repeat the arguments of~\refp{family-on-strata} (without dualization) .
In this way,  on $\URt^k$, we obtain an  exact sequence
\begin{equation*}
0 \to \CV'_\bu \to \CV_\bu \to \CW_\bu \to 0,
\end{equation*}
where $\CW_\bu$ is a $(k,2k,k)$-dimensional family of Artin representations, resp.
$\CV'_\bu$  a $\tnt$-semistable family of $\alpha(r,d,n-k)$-dimensional representations, of the quiver $\bpt$.
Applying the arguments of the proof of~\refl{stratification} to the family $\CV'_\bu$, we obtain a natural stratification
\begin{equation*}
\URt^{k,,\bu} = \URt^{k,\ge 0} \supset \URt^{k,\ge 1} \supset \URt^{k,\ge 2} \supset \dots \supset \URt^{k,\ge n-k} 
\supset \URt^{k,\ge n+1-k} = \emptyset,
\end{equation*}
by the length of the Artin sheaf $\Coker((\CV'_2)^\vee \boxtimes\CO \to (\CV'_1)^\vee \boxtimes\CO(1))$.
Furthermore,
on the stratum $\URt^{k,l} := \URt^{k,\ge l} \setminus \URt^{k,\ge l+1}$
we get  an exact sequence of representations
\begin{equation*}
0 \to \CW'_\bu \to \CV'_\bu \to \CU_\bu \to 0
\end{equation*}
where $\CW'_\bu$ is a $(l,2l,l)$-dimensional family of Artin representations, and
$\CU_\bu$ is a $\tnt$-semistable family of $\alpha(r,d,n-k-l)$-dimensional supermonadic representations.

For  $m=0,1,\ldots$, we define
\begin{equation*}
\URt^{\ge m} := \bigsqcup_{k + l \ge m} \URt^{k,l}.
\end{equation*}
It is clear that each $\URt^{\ge m}$ is a $\GL(\alpha(r,d,n))$-invariant  closed subset of $\URt$.
The corresponding closed subsets of  $\URt/\!/_{\tet^0}\GL(\alpha(r,d,n))$, 
the GIT quotient, give a stratification 
of the Uhlenbeck moduli space $\UMt(r,d,n)$.
Finally, we replace each stratum by its reduced underlying scheme.
\epr

%
%

\rem{another-stratification}
We could start with splitting of Artin subrepresentations first (and define in this way closed subsets $\URt^{\bullet,\ge l} \subset \URt$)
and then continue with splitting Artin quotient representations. Note that this will give {\em different}\/ two-index
stratification of the space $\URt$, but the resulting total stratification will be the same.
\erem

\rem{urt-grt}
We always have an embedding $\GRt^k \subset \URt^{0,k}$. Moreover, if $r$ and $d$ are coprime then this inclusion becomes an equality
\begin{equation*}
\GRt^k = \URt^{0,k}.
\end{equation*}
Indeed, the union of $\URt^{0,k}$ over all $k$ parameterizes all monadic $\tet^0$-semistable representations,
and by~\refc{coprime-stable} these are precisely all $\tnt$-semistable representations.
\erem

Before we go to the proof of the main theorem we need one more result.
We will use the notation of the proof of~\refl{stratification-uhlenbeck}.
Consider the stratum $\URt^{k,0}$ of the space $\URt$. We checked in the proof
of~\refl{stratification-uhlenbeck} that the universal representation over it
fits into an exact sequence
\begin{equation}\label{uvw-sequence}
0 \to \CU_\bu \to \CV_\bu \to \CW_\bu \to 0
\end{equation}
with $\CU_\bu$ being supermonadic and $\CW_\bu$ being Artin.

\lem{sequence-splits}
The subset $\URt^{k,\mathrm{split}} \subset \URt^k$ of all points over which the exact sequence~\eqref{uvw-sequence} splits is closed.
\elem
\prf
Indeed, $\URt^{k,\mathrm{split}} = \URt^k \cap \URt^{\ge k,\bullet} \cap \URt^{\bullet,\ge k}$
and both $\URt^{\ge k,\bullet}$ and $\URt^{\bullet,\ge k}$ are closed subsets in $\URt$ by their construction.
\epr

\ssec{}{Proof of~\reft{str}}

%
%
(1) Follows immediately from the inclusion $\GRt^k \subset \URt^{0,k}$.

(2) Equality $\UMt^0(r,d,n) = \GMt^0(r,d,n)$ for $r$ and $d$ coprime follows immediately
from the equality $\GRt^0 = \URt^{0,0}$ which is the only component of the stratum of $\URt$ giving $\UMt^0(r,d,n)$.

%
%

Further, consider the subset $\URt^{k,\mathrm{split}} \subset \URt^k$.
By~\refl{sequence-splits} it is closed.
It is also $\GL(\alpha(r,d,n))$-invariant. Finally, each point in $\URt^k$ is S-equivalent to a point in $\URt^{k,\mathrm{split}}$.
Therefore
%
%
%
\begin{equation*}
\UMt^k(r,d,n) = \URt^k/\!/_{\tet^0} \GL(\alpha(r,d,n)) = \URt^{k,\mathrm{split}}/\!/_{\tet^0} \GL(\alpha(r,d,n)).
\end{equation*}
So, it is enough to find a direct product decomposition for the right-hand side of the equality.
For this recall that the universal family of representations restricted to $\URt^{k,\mathrm{split}}$
splits canonically as a direct sum
\begin{equation*}
\CV_\bu = \CU_\bu \oplus \CW_\bu
\end{equation*}
of a supermonadic $\tnt$-semistable $\alpha(r,d,n-k)$-dimensional representation and of an Artin representation of dimension $(k,2k,k)$.
This decomposition gives a $\GL(\alpha(r,d,n))$-equivariant morphism to the homogeneous space
\begin{equation*}
\URt^{k,\mathrm{split}} \to \GL(\alpha(r,d,n)) / (\GL(\alpha(r,d,n-k) \times \GL(k,2k,k))
\end{equation*}
such that the fiber over a point is the product $\hsRt^0(r,d,n-k) \times \aRt(k,2k,k)$. It follows that
\begin{multline*}
\URt^{k,\mathrm{split}} /\!/_{\tet^0}  \GL(\alpha(r,d,n))
\\ \cong (\hsRt^0(r,d,n-k) \times \aRt(k,2k,k)) /\!/_{\tet^0}  (\GL(\alpha(r,d,n-k) \times \GL(k,2k,k))
\\ \cong (\hsRt^0(r,d,n-k) /\!/_{\tet^0}  \GL(\alpha(r,d,n-k)) \times (\aRt(k,2k,k)) /\!/_{\tet^0}  \GL(k,2k,k))
\\ \cong \UMt^0(r,d,n-k) \times \aMt(k,2k,k).
\end{multline*}
This together with~\refp{artin-moduli-k} proves part (2) of the Theorem.

(3) 
A split representation in the closure of the $\GL(\alpha(r,d,n))$-orbit of a $\tnt$-semistable representation $V_\bu(E)$
is isomorphic to a direct sum  $V'_\bu\oplus V''_\bu$, where $V'_\bu$ is the supermonadic quotient, resp.
$V''_\bu$ is the maximal Artin subrepresentation, of $V_\bu(E)$. By~\refp{tnt-ss},
the direct summands  $V'_\bu$ and $V''_\bu$ correspond to the sheaves $E^{**}$ and $E^{**}/E$, respectively.
The result follows.
\qed

\end{document}